\documentclass[11pt,reqno]{article}
\usepackage{amssymb, amsmath, amsthm, amsfonts, amscd, epsfig}
 \usepackage{color}
\usepackage[english]{babel}
\usepackage{bm}
\usepackage{graphicx, epsfig}
\usepackage{kotex,setspace}

\newtheorem{theorem}{Theorem}[section]
\newtheorem{cor}[theorem]{Corollary}
\newtheorem{lemma}[theorem]{Lemma}
\newtheorem{prop}[theorem]{Proposition}
\newtheorem{definition}{Definition}

\newtheorem{remark}{Remark}

\def\ep{\epsilon}

\newcommand{\Bp}{\mathbf{p}}

\newcommand{\Bx}{\mathbf{x}}

\newcommand{\RR}{\mathbb{R}}
\newcommand{\CC}{\mathbb{C}}

\newcommand{\p}{\partial}

\newcommand{\ds}{\displaystyle}
\newcommand{\eqnref}[1]{(\ref {#1})}
\newcommand{\pf}{\medskip \noindent {\sl Proof}. \ }
\renewcommand{\qed}{\hfill $\Box$ \medskip}
\newcommand{\Bu}{\mathbf{u}}

\newcommand{\beq}{\begin{equation}}
\newcommand{\eeq}{\end{equation}}
\newcommand{\RN}[1]{%
  \textup{\uppercase\expandafter{\romannumeral#1}}%
}

\numberwithin{equation}{section}
\numberwithin{figure}{section}

\begin{document}

\title{Stress concentration for two nearly touching circular holes  \thanks{\footnotesize
This work is supported by the Korean Ministry of Science, ICT and Future Planning through NRF grant Nos. 2013R1A1A3012931 (to M.L) and 2016R1A2B4014530 (to M.L).}}

\author{Mikyoung
Lim\thanks{\footnotesize Department of Mathematical Sciences,
Korea Advanced Institute of Science and Technology, Daejeon
305-701, Korea (mklim@kaist.ac.kr).} \and
Sanghyeon Yu\thanks{\footnotesize Seminar for Applied Mathematics, Department of Mathematics, ETH Z\"urich, R\"amistrasse 101, CH-8092 Z\"urich, Switzerland (sanghyeon.yu@sam.math.ethz.ch).}}

\maketitle

\begin{abstract}
We consider the plane elasticity problem for two circular holes. When two holes are close to touching, the stress concentration happens in the narrow gap region. In this paper, we characterize the stress singularity between the two holes by an explicit function.  
A new method of a singular asymptotic expansion for the Fourier series with slowing decaying coefficients is developed to investigate the asymptotic behavior of the stress.
\end{abstract}



\section{Introduction}

Stress concentration is a long studied subject in elasticity due to its practical importance. When two elastic inclusions are nearly touching, the stress distribution {can be arbitrarily large} in the narrow gap region. This blow-up {phenomenon} of the stress makes it challenging to numerically compute the field distribution. In this paper we analytically characterize the stress singularity {between close-to-touching} circular holes.

Let us discuss a similar problem in the context of conductivity. Consider two conducting inclusions which are separated by a distance $\ep$. We assume that the conductivity of the two inclusions is $k$.
Let $u$ denote the electric potential, which is a solution to the Laplace's equation, generated by the two inclusions.
As the two inclusions get closer, the resulting electric field $\nabla u$ can be very large in the gap region between the inclusions.
In fact, the asymptotic behavior of $\nabla u$ for small $\ep$ crucially depends on the conductivity $k$. 
If $k$ stays away from $0$ or $\infty$, then $|\nabla u|$ remains bounded regardless of how small $\ep$ is  (\cite{LV-ARMA-00}).
On the contrary, when $k=\infty$ or $k=0$, the electric field $|\nabla u|$ may blow up in the gap region as $\ep\rightarrow 0$.
Indeed, it was proved that in two dimensions, $|\nabla u|$ blows up as $\ep^{-1/2}$ when $k=\infty$ or $0$ (\cite{AKLLL-JMPA-07, AKLLZ-JDE-09, Yun-SIAP-07, Yun-JMAA-09, BLY-CPDE-10, LY-JDE-11, Gorb-MMS-16}). 
In three dimensions, the generic blow-up rate is $(\ep|\ln\ep|)^{-1}$ when $k=\infty$ 
(\cite{LY-CPDE-09, Lekner-PRSA-12, Gorb-MMS-16})
and $\ep^{-{(2-\sqrt 2)}/{2}}$ when $k=0$ (\cite{Yun-arXiv}).
It is worth to mention that a similar blow-up estimate was derived for the p-Laplace equation in \cite{GN-MMS-12}.

Now we return to the elasticity problem.
In contrast to the conductivity case, there are only a few results in the linear elasticity, {\it i.e.}, the Lam\'e system.
The difficulties come from both the vectorial nature of the elasticity and the fact that useful properties, such as the maximum principle, of a solution to Laplace's equation are not applicable to the Lam\'e system. 

Before considering the previous works on the elasticity, we introduce some definitions.
Let $B_1$ and $B_2$ be two disjoint elastic inclusions in $\mathbb{R}^2$ and let $\ep:=\mbox{dist}(B_1,B_2)$. 
We also assume that both the inclusions $B_1\cup B_2$ and the background are occupied by isotropic and homogeneous materials. Let $(\widetilde{\lambda},\widetilde{\mu})$ and $(\lambda,\mu)$ be the Lam\'{e} constants of the inclusions $B_1\cup B_2$ and of the background, respectively. 
For a displacement field $\mathbf{u}$, we define the stress tensor $\mathbf{\sigma}$ as
$${\sigma}:= \lambda (\nabla \cdot \mathbf{u})\mathbf{I} +\mu( \nabla \mathbf{u}+ (\nabla \mathbf{u})^T).
$$
Here, $\mathbf{I}$ is the $2\times 2$ identity matrix and the superscript $T$ denotes the transpose of a matrix.
It was proved in \cite{LN-CPAM-03} that, if the Lam\'e constants $\widetilde\lambda$ and $\widetilde\mu$ of the inclusion are finite, the stress $|\sigma|$ stays bounded regardless of $\ep$ (in fact, the result was proved for a more general elliptic system). 
But the situation is different when the Lam\'e constants $(\widetilde\lambda,\widetilde\mu)$ are extreme.
There are two extreme cases: 
hard inclusions $(\widetilde{\lambda}=\mbox{const.}, \widetilde{\mu}=\infty)$ and
holes $(\widetilde{\lambda}=0,\widetilde{\mu}=0)$.
For both cases, the stress $|\sigma|$ can be arbitrarily large in the gap region as $\ep\rightarrow 0$. 
The blow-up feature of the stress tensor was numerically verified in \cite{KK-CM-16}.

We now discuss previous works on the extreme cases.
For two hard inclusions which have general shapes in two dimensions, 
Li, Li and Bao \cite{BLL-ARMA-15} derived the $\ep^{-1/2}$ upper bound of the gradient $|\nabla \mathbf{u}|$. They also obtained the upper bounds for higher dimensions in \cite{BLL-AM-17}.
{
For two general-shaped holes, Bao, Li and Yin \cite{LY-CPDE-09} established the $\ep^{-1/2}$ upper bound.}
When the inclusions $B_1$ and $B_2$ are two circular holes in two dimensions, 
Callias and Markenscoff \cite{CM88,CM93} derived an asymptotic expansion of the stress $\sigma$ on the boundaries of the inclusions, by developing a singular asymptotic method. They also showed that the stress $\sigma$ blows  up as ${\ep}^{-1/2}$ as $\ep$ tends to zero. 
See also \cite{WM96}.

The purpose of this paper is to quantitatively characterize the stress singularity between two 2D circular holes under a uniform normal loading. 
Specifically, we derive an asymptotic expansion of the stress over the whole region outside the inlcusions.
As a result, we find an explicit function which completely captures the singular behavior of the stress distribution $\sigma$.
To our best knowledge, this is a first result of completely charaterizing the stress concentration for the hole case.

We shall see that the stress $\sigma$ is represented in terms of Fourier series with slowly {decaying} coefficients.  A new singular asymptotic expansion method is developed to deal with such series.
We emphasize that our method is completely different from the one in \cite{CM88,CM93,WM96}. 
Unfortunately, it seems that their method cannot be applied for a complete characterization of stress concentration (see Remark \ref{rmk:delicate_theta} in subsection \ref{sec:motivation}). We also emphasize that our approach is much simpler and provides important insights into an asymptotic behavior of the {Fourier} series.


{
It is worth to mention that, in \cite{MM}, an asymptotic solution for two circular elastic inclusions  was derived using a continuous distribution of point sources.
There, it was shown that high order multipole coefficients of their asymptotic solution are in good agreement with numerical results. However, their solution is not sufficiently accurate to capture the stress singularity in the gap region.
}



The paper is organized as follows. 
In section \ref{sec:main_result},  we state our main result. 
In section \ref{sec:Airy}, we review the Airy stress function formulation in the bipolar coordinates and then present an exact analytic solution for two circular holes derived by Ling \cite{CBLing}. In section \ref{sec:singular}, we propose a new method of singular asymptotic expansion. In section \ref{sec:stress}, we apply the proposed method for singular asymptotic expansion to Ling's analytic solution and then derive an asymptotic expansion of the stress tensor for two circular holes in the nearly touching limit. 


\section{Statement of results}\label{sec:main_result}

\begin{figure}[ht!]
\begin{center}
\epsfig{figure=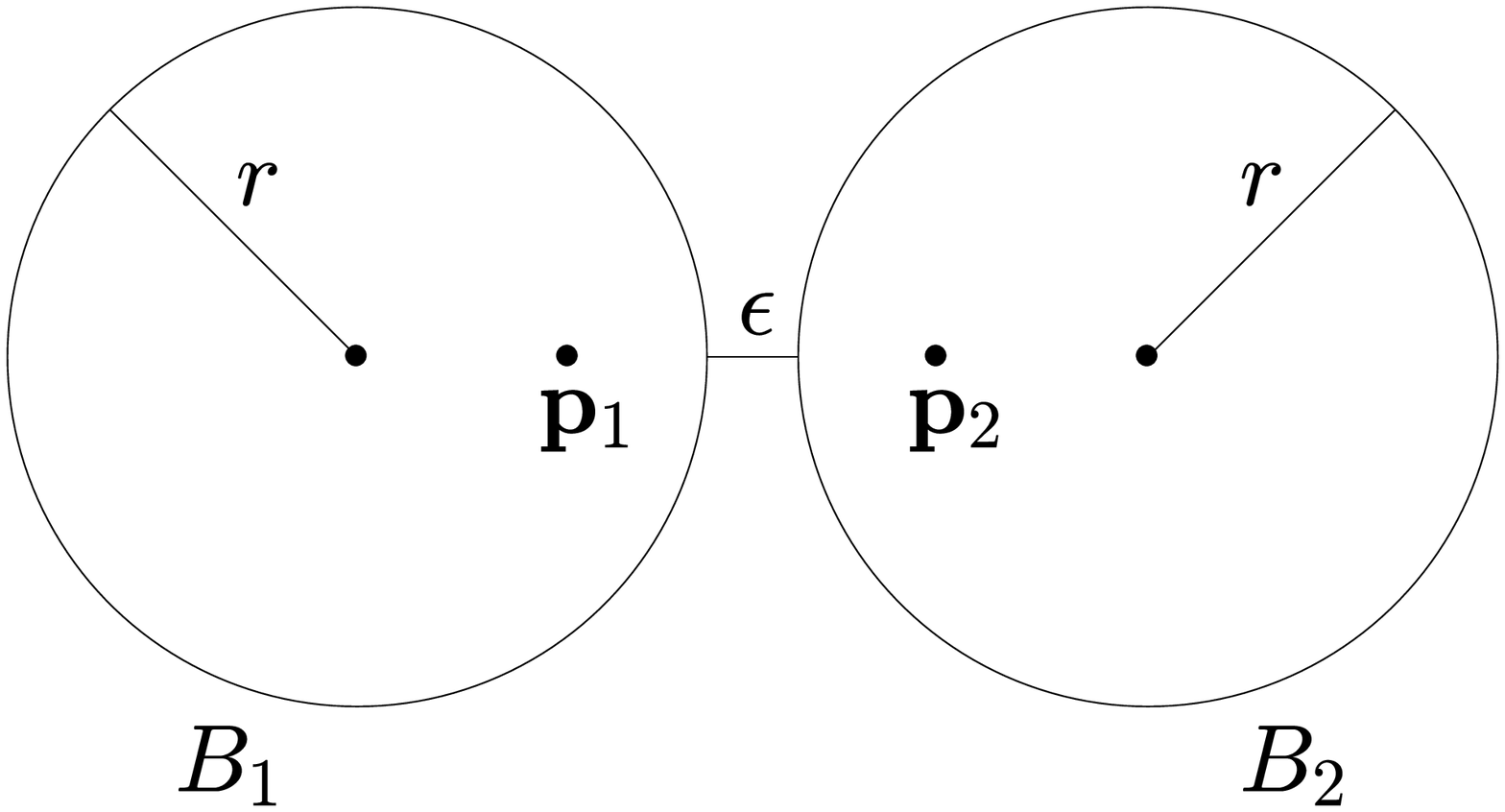,width=7.1cm}
\epsfig{figure=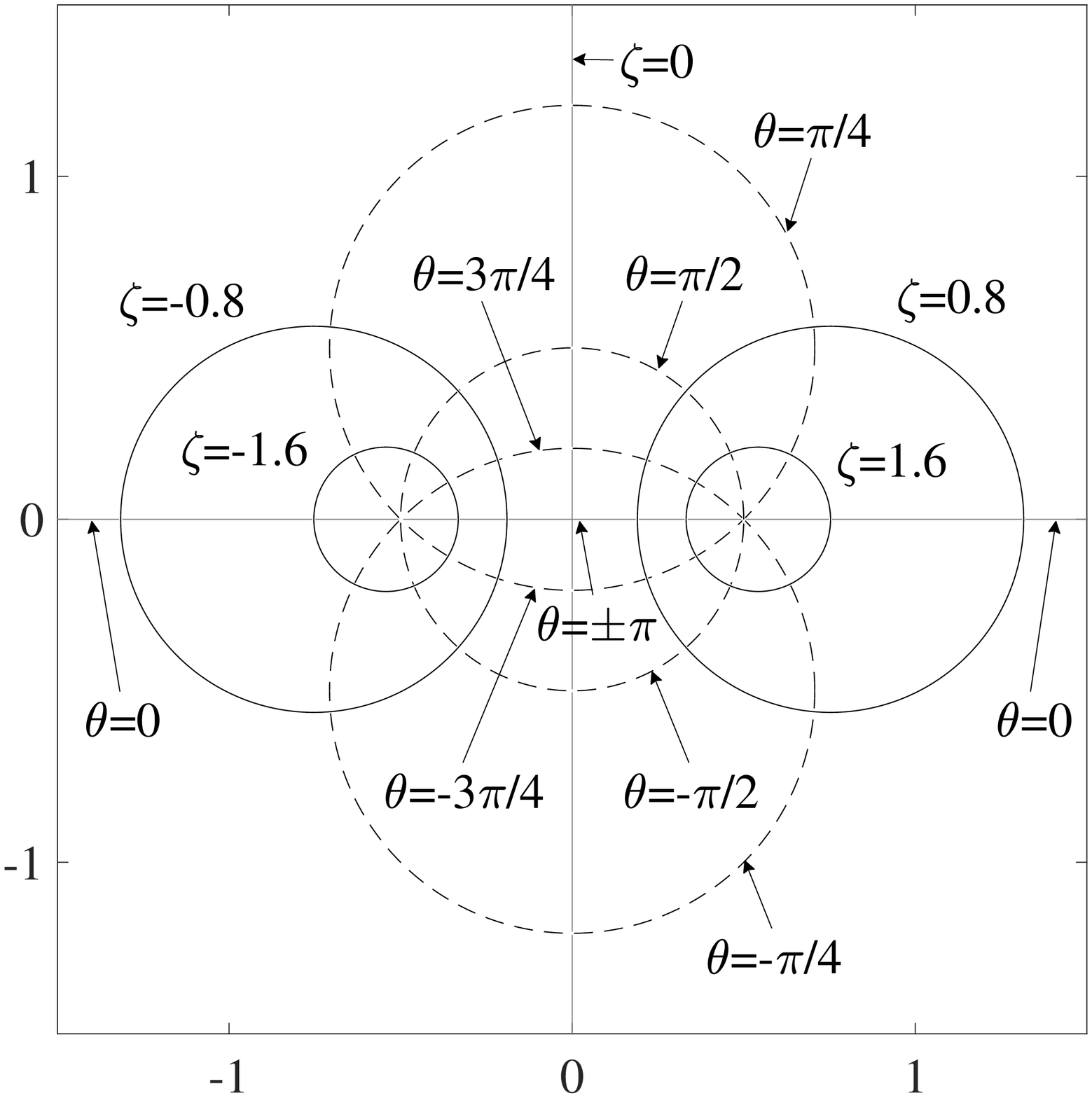,width=7.7cm}
\end{center}
\caption{(left) Geometry of the two circular holes, (right) The bipolar coordinates with $\alpha=0.5$}
\label{fig:bipolar}
\end{figure}

We assume that the inclusions $B_1$ and $B_2$ are circular disks of the same radius $r$. 
We also assume the inclusions are holes, {\it i.e.},
$
\widetilde\lambda=\widetilde\mu=0.
$
We may assume that the disks $B_1$ and $B_2$ are centered at $(-r-\ep/2,0)$ and $(r+\ep/2,0)$, respectively. {See Figure \ref{fig:bipolar}.}

The differential operator $\mathcal{L}_{\lambda,\mu}$ for the Lam\'{e} system is defined by
$$
\mathcal{L}_{\lambda,\mu} \mathbf{u}:= \mu \Delta \mathbf{u}+(\lambda+\mu)\nabla \nabla \cdot \mathbf{u}.
$$
Suppose that the Lam\'e constants $(\lambda,\mu)$ satisfy $\mu > 0$ and $\lambda+\mu>0$. Then $\mathcal{L}_{\lambda,\mu}$ becomes an elliptic operator.
The displacement field $\mathbf{u}$ is a solution to the Lam\'e system $\mathcal{L}_{\lambda,\mu} \mathbf{u}=0$ when the body force is absent. 
The conormal derivative  $\p \mathbf{u}/\p \nu$ on $\p B_j$ is given by
$$
\frac{\p\mathbf{u}}{\p \nu}:=\sigma  \nu = \lambda (\nabla \cdot \mathbf{u})\nu +\mu( \nabla \mathbf{u}+ (\nabla \mathbf{u})^T)\nu,$$
where $\nu$ is the outward unit normal vector to $\p B_j$.

When two circular holes $B_1 \cup B_2$ are embedded in the free space $\mathbb{R}^2$, the displacement field $\mathbf{u}$ satisfies the following equation:
\begin{equation} \label{elas_eqn}
 \ \left \{
 \begin{array} {ll}
\ds \mathcal{L}_{\lambda,\mu} \mathbf{u}= 0 \quad &\mbox{ in } \mathbb{R}^2 \setminus \overline{B_1 \cup B_2},\\[2mm]
\ds \frac{\p\mathbf{u}}{\p \nu}\Big|_+ =0 \quad &\mbox{ on } \p B_1 \cup \p B_2,\\[2mm]
\ds \mathbf{u}(\Bx)-\mathbf{u}_0(\Bx) = O(|\Bx|^{-1}) \quad& \mbox{ as } |\Bx| \rightarrow \infty,
 \end{array}
 \right.
 \end{equation}
where $\mathbf{u}_0$ is any solution to $\mathcal{L}_{\lambda,\mu} \mathbf{u}_0= 0 \mbox{ in } \mathbb{R}^2$ and the subscript $+$ denotes the limit from outside $\p B_j$. 
In this paper, we assume that $\Bu_0$ is a uniform normal loading given by
\beq\label{unif_normal_load}
\mathbf{u}_0(x,y)=\frac{1}{2(\lambda+\mu)}
\left[\begin{matrix}
x\\ y
\end{matrix}\right].
\eeq
One can easily check that the corresponding stress tensor $\sigma[\Bu_0]$ is the $2 \times 2$ identity matrix.

{
In this paper, we look for a decomposition of the stress tensor $\sigma=\sigma[\mathbf{u}]$ of the form
$$
\sigma = \sigma_* + \sigma_b, \quad \mbox{in } \mathbb{R}^2 \setminus (B_1 \cup B_2)
$$
such that $\sigma_*$ is an explicit function and $|\sigma_b|$ is bounded regardless of how small the distance $\ep$ is. Then we can say that $\sigma_*$ completely characterize the stress concentration. The goal is to find the function $\sigma_*$ explicitly.

}


{To state our result, we need some definitions.}
Let us define $\mathbf{p}_1$ and $\mathbf{p}_2$ as
\beq\label{def:p1p2}
\Bp_1 = (-\sqrt{\epsilon(r+\frac{\ep}{4})},\, 0\big) \quad \mbox{and}\quad
\Bp_2 = (\sqrt{\epsilon(r+\frac{\ep}{4})},\,0),
\eeq
and define a constant $\mathcal{I}_0$ as
\beq\label{def:I_0}
\mathcal{I}_0=\frac{1}{4}\int_{0}^\infty \frac{\sinh^2 x-x^2}{x^3 (\sinh 2x+ 2x)}dx.
\eeq
Let us denote $\Bx^{\perp} = (-y,x)$ for $\Bx=(x,y)$. Let $|\cdot|$ be the euclidean norm in $\mathbb{R}^2$. Let $\{\mathbf{e}_x,\mathbf{e}_y\}$ be the standard basis for $\mathbb{R}^2$.

The following is our main result in this paper {(for its proof, see subsection \ref{subsec:asymp_stress_tensor})}. 
The stress singularity between two nearly touching circular holes is explicitly characterized. 

\begin{theorem}\label{thm:main_stress_asymp}
Let $\Bu$ be the displacement field which is the solution to the elasticity problem \eqnref{elas_eqn} with \eqnref{unif_normal_load}.
Then its associated stress tensor $\sigma$ has the following decomposition:
\beq
\notag
\sigma = 
 \frac{r}{\mathcal{I}_0}|\mathbf{v}|(\frac{{\mathbf{w}}}{|{\mathbf{w}}|} \otimes \frac{{\mathbf{w}}}{|{\mathbf{w}}|}) + \sigma_b, 
\quad \mbox{ in }\mathbb{R}^2\setminus (B_1\cup B_2),
\eeq
where $\mathbf{v}$ and $\mathbf{w}$ are given by
\begin{align*}
\mathbf{v}(\Bx) &= \frac{\Bx-\Bp_1}{|\Bx-\Bp_1|^2} -\frac{\Bx-\Bp_2}{|\Bx-\Bp_2|^2}\quad \mbox{and}\quad
\mathbf{w}(\Bx) = \frac{(\Bx-\Bp_1)^\perp}{|\Bx-\Bp_1|^2} - \frac{(\Bx-\Bp_2)^\perp}{|\Bx-\Bp_2|^2}.
\end{align*}
Moreover, $\sigma_b$ satisfies
$$
\| \sigma_b \|_{L^\infty(\mathbb{R}^2\setminus(B_1\cup B_2))} \leq C,
$$
where $C>0$ is a constant independent of $\ep>0$.
\end{theorem}

\begin{cor}\label{cor_main_1}
{Under the same assumptions as in Theorem \ref{thm:main_stress_asymp},
the stress tensor $\sigma$ shows the following blow-up behavior at the origin}: for small $\ep>0$,
$$
\sigma(0,0) = \frac{2\sqrt{r}}{\mathcal{I}_0}\frac{1}{\sqrt{\ep}} \mathbf{e}_y \otimes \mathbf{e}_y + O(1).
$$
\end{cor}
{
\pf 
From the definitions of $\mathbf{p}_1,\mathbf{p}_2$ and $\mathbf{v}$, it is easy to check that $\mathbf{v}(0,0)=2/\sqrt{r\ep}\mathbf{e}_x + O(\ep^{3/2})$. Similarly, we have $(\mathbf{w}/|\mathbf{w}|)(0,0) = \mathbf{e}_y$. So the conclusion immediately follows from Theorem \ref{thm:main_stress_asymp}.
\qed
}

\begin{cor}
Under the same assumptions as in Theorem \ref{thm:main_stress_asymp},
the optimal blow-up rate of the stress tensor $\sigma$ is $\epsilon^{-1/2}$. More precisely, we have the following blow-up estimate: 
$$
\frac{C_1}{\sqrt\epsilon}\leq\| \sigma \|_{L^\infty(\mathbb{R}^2\setminus (B_1\cup B_2))} 
\leq \frac{C_2}{\sqrt\epsilon},
$$
for some positive  constants $C_1$ and $C_2$ independent of $\epsilon$.
\end{cor}
{\pf
Thw lower bound follows from Corollary \ref{cor_main_1}. 
It is easy to check that $|\mathbf{v}(x,y)| \leq C {\ep}/({\ep+y^2}) \leq C \ep^{-1/2}$ for $(x,y)\in\mathbb{R}^2\setminus(B_1\cup B_2)$. Here, $C$ is a positive constant independent of $\ep$. So we get the upper bound. The proof is completed.
\qed}

\begin{remark}
It is also important to consider the case of a shear loading $\Bu_0(x,y) \propto [y, x]^T$. 
While we only consider a uniform normal loading $\Bu_0(x,y) \propto [x,y]^T$ in this paper, our approach can be applied to the shear loading case as well. It will be investigated in a forthcoming paper.
\end{remark}

\section{Airy stress function for two circular holes}\label{sec:Airy}



\subsection{The bipolar coordinates}\label{section:bipolar}
We introduce the bipolar coordinates system and its properties. 
For a positive constant $\alpha>0$,
the bipolar coodinates $(\zeta,\theta)\in\RR\times (-\pi,\pi]$ is defined as
 \beq\label{xitheta1}
{\zeta-i\theta}=\log \frac{x+iy+\alpha}{x+iy-\alpha}.
\eeq
By separating \eqnref{xitheta1} into real part and imaginary part, it is easy to see that
\begin{equation}\label{xy}
x=\alpha\frac{ \sinh \zeta }{\cosh \zeta - \cos \theta}\qquad\mbox{and}\qquad y=\alpha\frac{ \sin \theta}{\cosh \zeta - \cos \theta}.
\end{equation}

The origin $(0,0)$ corresponds to $\zeta=0,\theta=\pm\pi$. The point at infinity corresponds to $(\zeta,\theta)=(0,0)$. It can be easily shown that the coordinate curves $\{\zeta=c\}$ and $\{\theta=c\}$ for a constant $c$ are respectively the zero level set of
\begin{align}
 f_\zeta(x,y,z)&=
\left(x-{\alpha}{\coth c}\right)^2 +y^2-\left({\alpha}/{\sinh c}\right)^2,
\label{fzeta}
\\
 f_\theta(x,y,z)&=
x^2+\left(y-{\alpha}{\cot c}\right)^2-\left({\alpha}/{\sin c}\right)^2.
\label{ftheta}
\end{align}
Note that the $\zeta$-coordinate curve is the circle of radius $\alpha/\sinh c$ centered at $(\alpha \coth c,0)$.
Therefore, $\zeta=c$ (resp. $\zeta=-c$) represents a circle contained in the region $x>0$ (resp. $x<0$). Moreover, $|\zeta|<c$ (resp. $|\zeta|>c$) represents the region outside (resp. inside) the two circles. 
Note also that $\theta$-coordinate curve is the circle of radius ${\alpha}/{\sin c}$ centered at $(\alpha \cot c,0)$. {See Figure \ref{fig:bipolar}.}

The boundaries $\p B_1$ and $\p B_2$  can be parametrized by $\{\zeta = -s\}$ and $\{\zeta = s\}$, respectively, for some suitable constant $s$.
Recall that $\p B_1$ and $\p B_2$ are the circles of the same radius $r$ centered at $(-r-\ep/2,0)$ and $(r+\ep/2,0)$, respectively.
In view of \eqnref{fzeta}, we choose $s$ and $\alpha$ such that $\alpha \coth s = r+\ep/2$ and $\alpha/\sinh s = r$.
Then one can easily check that
\beq\label{alpha_and_s}
\alpha= \sqrt{\epsilon(r+\frac{\ep}{4})} \quad \mbox{and}\quad s = \sinh^{-1} \sqrt{\frac{\epsilon}{r}(1+\frac{\ep}{4r})}.
\eeq
Note that $\alpha=O(\ep^{1/2})$ and $s=O(\ep^{1/2})$ for small $\ep>0$.
Now the boundaries $\p B_j$ and the regions $B_j$ can be represented as
\begin{align*}
&\p B_1 = \{ \zeta = -s\},
\quad
B_1 = \{ \zeta < -s\},
\\
&\p B_2 = \{ \zeta = +s\},
\quad
B_2 = \{ \zeta > +s\},
\end{align*}
and the exterior region outside $B_1\cup B_2$ is represented as
$$
\mathbb{R}^2\setminus (B_1 \cup B_2) = \{ |\zeta| \leq s\}.
$$


Let $\{\mathbf{e}_\zeta,\mathbf{e}_\theta\}$ be the unit basis vectors in the bipolar coordinates. Since the coordinate system is orthogonal, we have
\beq\label{unit_basis}
\mathbf{e}_\zeta = \nabla \zeta/|\nabla \zeta|, 
\quad
\mathbf{e}_\theta = \nabla \theta/|\nabla \theta|. 
\eeq
For a scalar function $f$, its gradient can be represented as
\beq\label{grad_f_bipolar}
\nabla f = \frac{\cosh \zeta-\cos\theta}{\alpha}\Big(\frac{\p f}{\p\zeta}\mathbf{e}_\zeta + \frac{\p f}{\p\theta}\mathbf{e}_\theta \Big).
\eeq
Here, $\nabla$ is the gradient with respect to $(x,y)$.

For later use, we define a function $h(\zeta,\theta)$ as
\beq \label{scale} 
h(\zeta,\theta):=\cosh\zeta-\cos\theta.
\eeq
The following two simple estimates regarding $h$ will be useful: for small $s>0$ and for all $|\theta|\leq \pi$,
\begin{align}
\label{s2_h_bdd}
&\frac{s}{[{h(s,\theta)}]^{1/2}} \leq \frac{s}{\sqrt{\cosh s-1}} \leq C,
\\
\label{sin_th_h_bdd}
&\frac{|\sin\theta|}{[h(s,\theta)]^{1/2}} \leq C\frac{\sqrt{1-\cos \theta}}{\sqrt{\cosh s-\cos\theta}} \leq C,
\end{align}
where $C$ is a positive constant independent of $s$.

\subsection{Airy stress function formulation}

We can reduce the equation \eqnref{elas_eqn} to a scalar problem.
It is well-known that, assuming that the body forces are negligible, the components of stress $\sigma$ can be represented as 
\beq\label{stress_cartesian}
\sigma_{xx}=\frac{\p^2 \chi}{\p y^2},\quad
\sigma_{xy}=\sigma_{yx}=-\frac{\p^2 \chi}{\p x\p y}, \quad
\sigma_{yy}=\frac{\p^2 \chi}{\p x^2}
\eeq
with a scalar function $\chi$  which satisfies the biharmonic equation $
\Delta \Delta \chi =0$. 
The function $\chi$ is called the {\it Airy stress function}. 

Let  $\chi_0$ be the stress function corresponding to the uniform normal loading $\mathbf{u}_0$. Then
it is easy to see that $
\chi_{0}(x,y) = \frac{1}{2}(x^2+y^2).
$
Let us decompose the total stress function $\chi$, which corresponds to $\Bu$, as follows:
$$
\chi = \chi_0 +  \chi_1.
$$
Then the stress function $\chi_1$ satisfies 
\begin{equation} \label{stress_eqn}
 \ \left \{
 \begin{array} {ll}
\ds \Delta \Delta \chi_1= 0 \quad &\mbox{ in } \mathbb{R}^2 \setminus \overline{B_1 \cup B_2},\\[1mm]
\ds \sigma[\chi_1+\chi_0]\nu|_+ =0 \quad &\mbox{ on } \p B_1 \cup \p B_2,\\[1mm]
\ds \sigma[\chi_1](\Bx)  \rightarrow 0 \quad &\mbox{ as } |\Bx| \rightarrow \infty,
 \end{array}
 \right.
 \end{equation}
 where $\sigma[\widetilde\chi]$
means the stress tensor associated to a stress function $\widetilde\chi$.

Let us denote
$$
\sigma^0 = \sigma[\chi_0]\quad\mbox{and}\quad
\sigma^1 = \sigma[\chi_1].
$$

In the following subsection, we shall present the exact analytic solution for the stress function $\chi_1$ and the corresponding stress tensor $\sigma^1$.

\subsection{Airy stress function in the bipolar coordinates}\label{subsec:AiryStress_bipolar}

In \cite{Jeff}, Jeffrey developed a general framework for plane elasticity problems in the bipolar coordinates. Let us briefly review their result. The biharmonic equation for the stress function $\chi$ can be written in terms of the bipolar coordinates as 
\beq\label{biharmonic_bipolar}
\Big(\frac{\p^4}{\p\theta^4}+2\frac{\p^4}{\p \zeta^2 \p\theta^2}+\frac{\p^4}{\p\zeta^4}+2\frac{\p^2}{\p \theta^2}-2\frac{\p^2}{\p \zeta^2}+1\Big) (\chi/J)=0
\eeq 
where $$J=\alpha/(\cosh\zeta-\cos\theta).$$
The components $\sigma_{\zeta\zeta},\sigma_{\theta\theta},\sigma_{\zeta\theta}$ of the stress tensor $\sigma$ in the bipolar coordinates are given by
 \beq\ds
\label{stress_formula}
\begin{cases}
\ds \sigma_{\zeta\zeta}[\chi]= \ds\frac{1}{\alpha}\Big[ (\cosh \zeta-\cos \theta)\frac{\p^2}{\p \theta^2} -\sinh{\zeta}\frac{\p}{\p\zeta} -\sin\theta \frac{\p}{\p\theta}+\cosh \zeta\Big]  (\chi/J),\\[1mm]
\ds \sigma_{\theta\theta}[\chi] = \ds\frac{1}{\alpha}\Big[ (\cosh \zeta-\cos \theta)\frac{\p^2}{\p \zeta^2} -\sinh{\zeta}\frac{\p}{\p\zeta} -\sin\theta \frac{\p}{\p\theta}+\cos \theta\Big] (\chi/J),\\[1mm]
\ds\sigma_{\zeta\theta}[\chi] = \ds\sigma_{\theta\zeta}[\chi]= -\frac{1}{\alpha}\Bigr[{(\cosh\zeta-\cos \theta)} \frac{\p^2}{\p \zeta \p \theta}\Bigr](\chi/J).
\end{cases}\eeq

Let us now consider the stress function $\chi_1$ for the two circular holes problem \eqnref{stress_eqn}. 
Its analytic expression was derived by Ling \cite{CBLing} as follows:
\begin{align}\label{Psi_a_series}
\ds \frac{1}{\alpha} (\chi_1/J )(\zeta,\theta)&=  K (\cosh \zeta-\cos \theta)\log( \cosh \zeta-\cos \theta) +\sum_{n=1}^{\infty}
\phi_n(\zeta) \cos(n \theta),
\end{align}
where
\beq\label{phin}
\phi_n(\zeta)=A_n\cosh(n+1)\zeta + B_n\cosh(n-1)\zeta.
\eeq
Here, $A_n$ and $B_n$  are constant coefficients given by
\beq\label{eqn:ABK}\ds
\begin{cases}
\ds
A_n =\ds\frac{2K (e^{-ns}\sinh ns+ne^{-s} \sinh s)
}{n(n+1)(\sinh 2ns+n\sinh 2s)}, \quad n\geq 1,
\\[3mm]
\ds B_n = \ds-\frac{2K (e^{-ns}\sinh ns+ne^{s} \sinh s)
}{n(n-1)(\sinh 2ns+n\sinh 2s)},\quad n \geq 2,\\[3mm]
\ds B_1=\ds\frac{1}{2}\left(K \tanh s \cosh 2s-2e^{-s} \cosh s
\right).
\end{cases}
\eeq
Also, the constant $K$ is given by
\begin{align}
K= \big(\frac{1}{2}+\tanh{s}\sinh^2{s}-  4 P(s) \big)^{-1}, \label{def_K}
\end{align}
where \beq\label{def_P}
P(s)=
\sum_{n=2}^{\infty}\frac{e^{-n s}\sinh ns+n(n\sinh s+\cosh s)\sinh s}{n(n^2-1)(\sinh 2ns + n\sinh 2s )}.
\eeq

Using the exact solution \eqnref{Psi_a_series} for $\chi_1$ and the stress components formulas \eqnref{stress_formula}, it is possible to derive explicit formulas for all components of $\sigma^1$.
We compute
\begin{align}
&(\sigma^1_{\theta\theta}-\sigma^1_{\zeta\zeta})(\zeta,\theta)\notag\\\notag
&=\frac{1}{\alpha}(\cosh\zeta-\cos\theta)\Big[ \frac{\p^2}{\p \zeta^2}-\frac{\p^2}{\p \theta^2}-1\Big](\chi_1/J )\\
\ds&=K\big(\cosh 2\zeta-2\cosh\zeta \cos\theta+\cos 2\theta\big)\notag\\\label{eqn:sigma1}
\ds&+h(\zeta,\theta)\sum_{n=1}^{\infty} \Big[2n(n+1)A_n\cosh(n+1)\zeta+2n(n-1)B_n\cosh(n-1)\zeta\Big] \cos n\theta.
\end{align}
Similarly, we get
\begin{align}
\sigma_{\zeta\theta}^1(\zeta,\theta)\notag
&
= 
-K\sinh\zeta\sin\theta 
\\&\quad 
+ h(\zeta,\theta)
\sum_{n=1}^\infty\Big[n(n+1)A_n\sinh(n+1)\zeta+n(n-1)B_n\sinh(n-1)\zeta\Big]\sin n\theta,\label{eqn:sigma2}
\\
\ds\sigma_{\zeta\zeta}^1(\zeta,\theta) &=\notag
-\frac{K}{2}(\cosh2\zeta-2\cosh \zeta \cos\theta +\cos 2\theta)
+
\sum_{n=1}^\infty\Big[-{h(\zeta,\theta)} n^2\phi_n(\zeta) \cos n\theta
 \\
 &\quad
 -\sinh \zeta\, \phi_n'(\zeta)  \cos n\theta
 +n \phi_n(\zeta)\sin \theta  \sin n\theta + \cosh\zeta \,\phi_n(\zeta) \cos n\theta)\Big].\label{eqn:sigma3}
\end{align}
The difference $(\sigma_{\theta\theta}^1-\sigma_{\zeta\zeta}^1)$ is considered instead of $\sigma_{\theta\theta}^1$ because it has a simpler expression.
In the next sections, we shall investigate the asymptotic behavior of the above series when the distance $\epsilon$ is small.

{The stress components admit much simpler expressions on the boundary $\p B_1 \cup \p B_2$.}
From the zero-traction condition, {\it i.e.}, {$\sigma\nu=0$} on the hole boundaries, we get $$\sigma_{\zeta\theta}^1=-\sigma_{\zeta\theta}^0\quad\mbox{and}\quad\sigma_{\zeta\zeta}^1=-\sigma_{\zeta\zeta}^0, \quad\mbox{on }\p B_1 \cup \p B_2.$$
For the component $\sigma_{\theta\theta}^1$, it was shown in \cite{CBLing} that
\begin{align}
\ds\sigma_{\theta\theta}^1\big|_{\p B_i} &= \sigma^1_{\theta\theta}(\zeta=(-1)^i s,\theta)
\notag
\\[.5mm]
\ds&=(2K \sinh s)\, h(s,\theta) \big[1+4q(s,\theta)\big],\quad i=1,2,
\label{stress:CBLing}
\end{align}
with 
\beq\label{sum:CBLing}
\ds q(s,\theta):=\sum_{n=1}^\infty \frac{\sinh n s}{\sinh 2n s + n \sinh 2 s}\cos n\theta.
\eeq

{We will need an asymptotic expansion of the constant $K$ {for small $s$}. }
We have the following lemma, whose proof is given in Appendix \ref{appendix:Ks}.
\begin{lemma}\label{lem:asymp_Ks}
For small $s>0$, we have
$$
K= \frac{1}{\mathcal{I}_0 s^2} (1+ O(s)),$$
where $\mathcal{I}_0$ is given by \eqnref{def:I_0}.

\end{lemma}



\section{Singular asymptotic expansion}\label{sec:singular}

In this section, we propose a new method of singular asymptotic expansion for infinite series.
We first explain, in section \ref{sec:motivation}, the motivation and main idea of our method by considering its simplified version. 
{Then, in section \ref{sec:singular_asymp}, we present a complete version of our method and its proof.}

 \subsection{Motivation and main idea}\label{sec:motivation}

{
Here we consider, for the ease of presentation, the stress component $\sigma_{\theta\theta}^1$ only on the boundary $\p B_2$(or $\{\zeta=s\}$). 
In the later sections, we will consider the stress tensor over the whole exterior region $\mathbb{R}^2\setminus{(B_1\cup B_2)}$. 
}

Recall that the analytic expression \eqnref{stress:CBLing} of $\sigma^1_{\theta\theta}(\zeta=s,\theta)$ contains the Fourier cosine series $q(s,\theta)$ given by
\beq
\ds q(s,\theta)=\sum_{n=1}^\infty q_n(s)\cos n\theta,\qquad\mbox{with}\quad q_n(s)=\frac{\sinh n s}{\sinh 2n s + n \sinh 2 s}.
\eeq
We are interested in the asymptotic behavior of the series $q(s,\theta)$ when the gap distance $\ep$ tends to zero. 
We shall use $s$ as a small parameter because $s$ is small as $s=O(\ep^{1/2})$.  

Throughout this paper, $C$ denotes a positive constant independent of $s>0$. 


\smallskip\smallskip
{\it\noindent\textbf{Difficulties in the nearly touching case.}}
Let us discuss difficulties in studying the series $q(s,\theta)$ when $s$ goes to zero.
A standard way to get {an asymptotic expansion of series} is to use the Taylor expansion.
Since 
\begin{align*}
q_n(s) &= q_n(0) + s q_n'(0) + \frac{s^2}{2} q_n''(0)+\cdots 
\\
&=\frac{1}{4} + s\cdot 0 - s^2 \frac{n^2+2}{24} +\cdots,
\end{align*}
we get the (formal) asymptotic formula
$$
\ds q(s,\theta)=\sum_{n=1}^\infty q_n(s)\cos n\theta = \sum_{n=1}^\infty \frac{1}{4} \cos n\theta -  s^2\sum_{n=1}^\infty \frac{n^2+2}{24} \cos n\theta + \cdots.
$$
Clearly, each term in the right hand side is {\it not} convergent. So, the above formal {formula} of $q$ fails to describe the asymptotic behavior. This originates from the slowly decaying property of $q_n(s)$:
\begin{center} $q_n(s)$ decays like $e^{-ns}$ as $n\rightarrow \infty$.
\end{center}
 In the limit $s\rightarrow 0$, the sequence $q_n(s)$ does not decay to zero, which results in non-convergence of the series in the above formula.


Besides the non-decaying feature of $q_n(s)$, the oscillating property of the cosine function makes even more difficult to understand the behavior for the series $q(s,\theta)$. 
{In fact, the asymptotic behavior of $q(s,\theta)$ for small $s$ can be dramatically different if the angle $\theta$ changes.} Note that $q(s,\theta)$ is an alternating series at $\theta=\pm\pi$, while it is a positive series at $\theta=0$. 
So we have 
\begin{align}
\ds &q(s,0) \approx \sum_{n=1}^\infty  \frac{1}{s}   \frac{s\cdot\sinh n s}{\sinh 2n s + 2n  s} 
\approx \frac{1}{s}\int_0^\infty \frac{\sinh x}{\sinh 2x+2x} dx \approx C\frac{1}{s},\label{Szero}\\
\ds
&q(s,\pm \pi) \approx \sum_{n=1}^\infty \frac{\sinh n s}{\sinh 2n s + 2n  s} (-1)^n 
\approx \frac{1}{2 s} \left( -\int_{s}^{2s} \frac{\sinh (x/2)}{\sinh x+ x} \,dx\right)
\approx -\frac{1}{8}.
 \label{Spmpi}
\end{align}
While $q(s,0)$ is as large as $s^{-1}$, $q(s,\pm\pi)$ converges to some constant as $s$ tends to zero.

For general {angles $\theta$}, Callias and Markenscoff proved in \cite{CM88, CM93} the following result using their singular asymptotic method for integrals:
\beq\label{CM_series}
q(s, \theta) = -\frac{1}{8}+O(s) \qquad\mbox{for fixed } \theta\neq0.
\eeq
However, Eq.\;\eqnref{CM_series} was obtained under the assumption that $\theta$ is a nonzero fixed constant and, therefore, it may not hold uniformly on $\{0<|\theta|\leq\pi\}$. In fact, Eq.\;\eqnref{CM_series} does not explain the transition of the asymptotic behavior of $q(s,\theta)$ from $O(1)$ to $ O(s^{-1})$ as $\theta$ tends to $0$.

We emphasize that the uniformity on $\theta$ is essential for understanding the stress concentration on the {\textquoteleft whole\textquoteright} boundary or in the {\textquoteleft whole\textquoteright} exterior domain, see Remark \ref{rmk:delicate_theta}. 


\smallskip
\smallskip
{\it\noindent{\textbf{Main idea of our approach}} }
Now we illustrate our approach to overcome the aforementioned difficulties. 
We first rewrite the coefficient $q_n(s)$ as a function of $ns$. Specifically, we write 
$$q(s,\theta)=\sum_{n=1}^\infty f(ns)\cos n\theta$$ with a smooth function $f$ given by
$$
f(x)=\frac{\sinh x }{\sinh 2x + c\, 2x},\qquad c=\frac{{\sinh 2s}}{{2s}}=1+O(s^2).$$ 
Note that $f(ns)$ decays as $e^{-ns}$ as $x\rightarrow \infty$. As a result, the standard approach fails because $\lim_{s\rightarrow 0} f(ns)$ does not decay to zero as $n\rightarrow \infty$ as already explained. 
Roughly speaking, our strategy for overcoming this difficulty is to consider the Taylor expansion of $f(x)e^{x}$ but not of $f(x)$.

Let us denote $$\widetilde{f}(x):=f(x)e^{x}$$ and set
 $$
 f_0:=\lim_{x\rightarrow 0+} {f}(x)e^{x}=\frac{1}{4} + O(s^2).
 $$
Then, $\widetilde{f}(x)$ has the following zeroth-order Taylor expansion at $x=0$:
$$
\widetilde{f}(x) = f_0 + r(x),
$$
where the remainder term $r(x)$ satisfies
$$
|r(x)|\leq \sup_{x_*\in (0,x)}|\widetilde{f}'(x_*)|\,x \leq 2x, \quad \mbox{for small } x>0.
$$ 
Therefore, the series $q(s,\theta)$ can be decomposed as
\begin{align*}
q(s,\theta)&
=\sum_{n=1}^\infty [f(ns)e^{ns} ]e^{-ns}\cos n\theta = \sum_{n=1}^\infty \widetilde{f}(ns)e^{-ns} \cos n\theta
\\
&=\sum_{n=1}^\infty f_0 e^{-ns}\cos n\theta +
\sum_{n=1}^\infty  r(ns)e^{-ns}\cos n\theta 
\\
&: = q_0(s,\theta) + q_1(s,\theta).
\end{align*}
Note that the series $q_0(s,\theta)$ and $q_1(s,\theta)$ are convergent contrary to the standard approach.
The leading order term $q_0(s,\theta)$ can be evaluated analytically to give
\beq\label{esti:q0}
q_0(s,\theta)=f_0\frac{(-1)}{2}\frac{e^{-s}-\cos\theta}{\cosh s-\cos\theta}.
\eeq 
We then need to consider the remainder term $q_1(s,\theta)$. In fact, it is tricky to derive an estimate of $q_1(s,\theta)$ because of its delicate dependence on $\theta$. 
Moreover, we should verify that $|q_1(s,\theta)|$ is smaller than $|q(s,\theta)|$ in a certain sense.

In order to gain some insight, it is better to estimate $q(s,\theta)$ before considering $q_1(s,\theta)$.
As already explained, the asymptotic behavior of $q(s,\theta)$ crucially depends on $\theta$. 
We need a good enough estimate for $q(s,\theta)$ so that it implies $q(s,0) = O(s^{-1})$ and $q(s,\pm\pi)=O(1)$.
Let us try to estimate ${q}(s,\theta)$ directly from its definition. We get
$$
|{q}(s,\theta) | \leq \sum_{n=1}^\infty |f(ns)| \leq C\frac{1}{s} \int_0^\infty |f(x)|dx \leq C\frac{1}{s}.
$$
Unfortunately, this estimate does not show the dependence of $q(s,\theta)$ on $\theta$. For example, it does not imply $q(s,\pm\pi)= O(1)$.

There is a simple but powerful way to overcome this difficulty. Let us consider a complex-valued version $\tilde{q}(s,\theta)$ given by 
$$\tilde{q}(s,\theta):=\sum_{n=1}^\infty f(ns)e^{in\theta}=\sum_{n=1}^\infty \widetilde{f}(ns)e^{n(-s+i\theta)}.$$
Note that $q(s,\theta)=\mbox{Re}\{\tilde{q}(s,\theta)\}$.
We consider $(1-e^{-s+i\theta}) \tilde{q}(s,\theta)$ instead of $\tilde{q}(s,\theta)$ and then rewrite it as
\begin{align*}
(1-e^{-s+i\theta}) \tilde{q}(s,\theta) &= \sum_{n=1}^\infty \left[\widetilde{f}(n s) e^{n(-s + i\theta)} - \widetilde{f}(ns) e^{(n+1)(-s + i\theta)}\right]
\\
&=\widetilde{f}(s)e^{-s+i\theta} + \sum_{n=2}^\infty \big[ \widetilde{f}(ns)-\widetilde{f}((n-1)s) \big] e^{-ns + in\theta}.
\end{align*}
Note that the second term in the RHS contains the difference $ \widetilde{f}(ns)-\widetilde{f}((n-1)s)$, which is smaller than $ \widetilde{f}(ns)$. So one can expect that a finer estimate can be obtained. Indeed, by the mean value theorem, we have 
$$
{ \widetilde{f}(ns)-\widetilde{f}((n-1)s)} = s \widetilde{f}'(s_n^*),\quad\mbox{for some } s_n^*\in ((n-1)s,ns).
$$
Therefore, we obtain
\begin{align*}
|(1-e^{-s+i\theta}) \tilde{q}(s,\theta)| & \leq |\widetilde{f}(s)|e^{-s} + \sum_{n=2}^\infty s |\widetilde{f}'(s_n^*)|e^{-ns}
\\
&\leq |{f}(s)| + C\sum_{n=2}^\infty s |f'(ns) e^{ns} + f(ns)  e^{ns})| e^{-n s}\\
&\leq |{f}(s)| + C \int_0^\infty (|f'(x)| + |f(x)|)dx \leq C.
\end{align*}
Since $|1-e^{-s+i\theta}|^2=2e^{-s}(\cosh s-\cos\theta)$, we get 
$$
|q(s,\theta)| \leq |\tilde{q}(s,\theta)|\leq C \frac{1}{\sqrt{\cosh s-\cos\theta}}.
$$
This estimate clearly shows the dependence of $q(s,\theta)$ on $\theta$.
Moreover, it implies both $q(s,0) = O(s^{-1})$ and $q(s,\pm\pi)=O(1)$, as desired.

Next we return to the remainder term $q_1(s,\theta)$. Let $\tilde{q}_1(s,\theta) $ be its complex version, namely, $\tilde{q}_1(s,\theta) = \sum_{n=1}^\infty r(ns)e^{-ns+in\theta}$ and we consider $(1-e^{-s+i\theta})^2 \tilde{q}_1(s,\theta)$. 
{Although we omit the details,}
it turns out that
a process similar to the case of $q(s,\theta)$ yields 
\begin{align*}
|(1-e^{-s+i\theta})^2 \tilde{q}_1(s,\theta)| 
\leq C s. 
\end{align*}
Hence, we obtain
$$
|q_1(s,\theta)| \leq C\frac{s}{\cosh s-\cos\theta}.
$$
Therefore, together with \eqnref{esti:q0}, we finally obtain an asymptotic expansion of the series $q(s,\theta)$ as follows: for small $s>0$,
\begin{align}
q(s,\theta)& = f_0\frac{(-1)}{2}\frac{e^{-s}-\cos\theta}{\cosh s-\cos\theta} + O(\frac{s}{\cosh s-\cos\theta})\notag\\
&=-\frac{1}{8}\frac{e^{-s}-\cos\theta}{\cosh s-\cos\theta} + O(\frac{s}{\cosh s-\cos\theta}).\label{asymp:q_onboundary}
\end{align}
{Note that this asymptotic formula holds uniformly on $\{|\theta|\leq\pi\}$ and it captures the transition of the asymptotic behavior of $q(s,\theta)$ from $O(1)$ to $ O(s^{-1})$ as $\theta$ tends to $0$.}
Consequently, from \eqnref{stress:CBLing}, an asymptotic expansion of the stress component $\sigma^1_{\theta\theta}(s,\theta)$ on the boundary $\p B_2$ immediately follows.

{
To summarize, we have shown how to derive an asymptotic expansion of the Fourier cosine series $q(s,\theta)$ for small $s$. 
In the next subsection, we will develop an asymptotic method for a general class of Fourier series with slowly decaying coefficients. It will enables us to investigate the stress tensor over the whole exterior region $\mathbb{R}^2\setminus (B_1\cup B_2)$.
}

\begin{figure*}
\begin{center}
\epsfig{figure=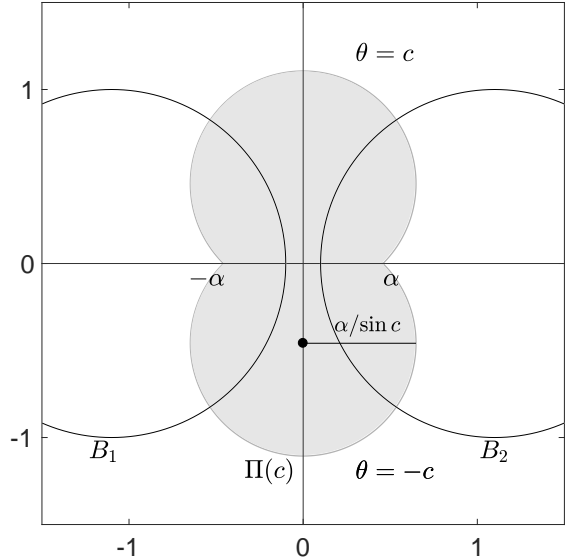,width=9.0cm}
\end{center}
\caption{Geometry of the region $\Pi(c)=\{c \leq|\theta|\leq \pi\}$. We set $r=1$ and $\epsilon=0.2$. The shaded region is $\Pi(c)$ when $c=\pi/4$.}
\label{fig:snowman}
\end{figure*}

\begin{remark}\label{rmk:delicate_theta}
We emphasize that the uniformity of the asymptotic expansion w.r.t. $\theta$ is important for the characterization of the stress concentration over the whole boundary $\p B_1 \cup \p B_2$ or over the whole exterior region $\mathbb{R}^2\setminus (B_1\cup B_2)$. 
{To see this, let us consider a region $\Pi(c):=\{(x,y):c \leq |\theta|\leq \pi\}$ for a constant $0 \leq c < \pi$. 
As can be seen in Figure \ref{fig:snowman}, the region $\Pi(c)$ consists of two intersecting disks of radius $\alpha/\sin c$, which contain the gap region betweeen the inclusions $B_1$ and $B_2$. The region $\Pi(c)$ has a delicate behavior when $\ep$ is small. Recall that $\alpha\approx\sqrt{r \ep}$. If we fix $c>0$ independently of $\ep$, then the radius $\alpha/\sin c$ is of order $ O(\ep^{1/2})$ and, thus, the region $\Pi(c)$ becomes vanishingly small as $\ep\rightarrow 0$. 
But, for a complete characterization of the stress over the whole gap region, the size of region of interest should remain essentially unchanged even when $\ep\rightarrow 0$.
In fact, we  should choose $c$ small enough to fix the size of $\Pi(c)$. 
For example, consider the case when $c=\sqrt\epsilon$.
In this case, the size of the region $\Pi(\sqrt{\ep})= \{ (x,y):\sqrt{\ep} \leq |\theta| \leq \pi \} $ does not vanish when $\ep$ goes to $0$ because $\alpha/\sin c \approx \sqrt r$.}
For this reason, it is needed to consider the case when $\theta$ is close to zero as well as $\pm\pi$.

\end{remark}

\subsection{New method for singular asymptotic expansion} \label{sec:singular_asymp}

Here we present our new method for singular asymptotic expansion {of
the Fourier series with slowly decaying coefficients.}

Let us denote $\mathbb{R}^+:=(0,\infty)$.
For a small parameter $s>0$ and a smooth function $f:\mathbb{R}^+\times (0,1) \rightarrow \mathbb{R}$, we define the complex function $\mathcal{L}[f]$ as
\beq\label{def:L_fs}
\mathcal{L}[f](z): = \sum_{n=1}^\infty f(n s,s) z^n,\qquad  z\in \mathbb{C}.
\eeq
{The series in the RHS is convergent for $|z|\leq e^s$ if $f$ satisfies $|f(x,y)|\leq C x^N e^{-2x}$ for all $1\leq x<\infty$, $y\in(0,1)$ and for some $N\in\mathbb{N}$.}

The following theorem is the main result in this section.
\begin{theorem} \label{thm:singular_L}
Let $f:\mathbb{R}^+\times (0,1)\rightarrow \mathbb{R}$ be a smooth function satisfying the following conditions:
\begin{itemize}
\item[\rm(A1)] For each $y\in(0,1)$, the limit $f_0(y):=\lim_{x\rightarrow 0+}f(x,y)$ exists.
\item[\rm(A2)] There exist a constant $C>0$ and a positive integer $N\in\mathbb{N}$ such that, for all $y\in(0,1)$ and $k=0,1,2$, it holds that
\beq\label{decay_fs}
\left|\p_x^{(k)} f(x,y)\right| \leq 
\begin{cases}
 C,&\quad 0<x<1,
\\[2mm]
C x^N e^{-2x},&\quad 1\leq x< \infty.
\end{cases}
\eeq 
\end{itemize}
Then, for small $s>0$, the series $\mathcal{L}[f](z)$ has the following asymptotic expansion:
$$
\mathcal{L}[f](z) = f_0(s)\frac{e^{-2s}z}{1-e^{-2s}z}  + O(\frac{s}{|1-e^{-2s}z|^2}), \qquad  |z|\leq e^{s}, \, z\in \mathbb{C}.
$$

\end{theorem}
\pf 
Fix $s$ to be a small positive number and let $\widetilde{f}(x):=f(x,s)e^{2x}$. From the assumption (A1), we have
$$
\lim_{x\rightarrow 0+}\widetilde{f}(x)=\lim_{x\rightarrow 0+}f(x,s)e^{2x}=f_0(s).
$$
We also have, from the assumption \eqnref{decay_fs}, that
$$
|{\widetilde{f}}'(x)| \leq \left|(\p_x f)(x,s) e^{2x} + f(x,s) 2 e^{2x} \right| \leq C, \quad \mbox{for } 0<x<1.
$$
Remind that $C$ denotes a positive constatnt independent of $s$.
So the function $\widetilde{f}(x)$ has the zeroth order Taylor expansion about $x=0$ as follows:
\beq\label{g_taylor}
\widetilde{f}(x) = f_0(s) + r(x),
\eeq
where the remainder term $r$ satisfies
\beq\label{rs_near_zero}
|r(x)| \leq \sup_{x_*\in (0,x)}|\widetilde{f}'(x_*)| \leq C x, \quad \mbox{for small } x>0.
\eeq
Note that
\beq\label{rdd_gdd_equal}
\widetilde{f}''(x) = r''(x).
\eeq

We now decompose $\mathcal{L}[f](z)$ using the Taylor expansion of $\widetilde{f}$. 
For $z\in \CC$ satisfying $|z|\leq e^s$, we obtain from \eqnref{g_taylor} that
\begin{align}
\mathcal{L}[f](z)
&=\sum_{n=1}^\infty \Big[f(ns,s) e^{2ns}\Big] e^{-2ns}z^n
=\sum_{n=1}^\infty \widetilde{f}(ns) e^{-2ns}z^n
\notag
\\
&=\sum_{n=1}^\infty f_0(s) e^{-2ns }z^n
+
\sum_{n=1}^\infty r(ns) e^{-2ns} z^n
\notag
\\
&=f_0(s)\frac{e^{-2s}z}{1-e^{-2s}z}
+
\sum_{n=1}^\infty r(ns)e^{-2ns}  z^n.
\end{align}

Now it only remains to show
\beq\label{def:remainderR}
R(z):=\sum_{n=1}^\infty r(ns) e^{-2ns} z^n = O(\frac{s}{|1-e^{-2s}z|^2}).
\eeq
For notational simplicity, let us denote $\tilde{z} = e^{-2s}z$ and $r_n = r(ns)$. To prove \eqnref{def:remainderR}, we consider $(1-\tilde z)^2 R$. We have
\begin{align*}
(1-\tilde z)^2 R(z)&= 
(1- 2\tilde{z} + \tilde{z}^2 ) \sum_{n=1}^\infty r_n \tilde{z}^n
\\&
= \sum_{n=1}^\infty r_n (\tilde{z}^n- 2 \tilde{z}^{n+1} + \tilde{z}^{n+2})
\\
&= \Big[r_1 \tilde{z} +(-2 r_1+r_2) \tilde{z}^2 \Big]
 + \sum_{n=3}^\infty \big(r_n -2 r_{n-1} +r_{n-2}\big) \tilde{z}^n : = I + II.
\end{align*}
We estimate $I$ and $II$ separately. From \eqnref{rs_near_zero} and the fact $|{z}|\leq e^{s}$, we get 
\begin{align*}
|I| &= |r(s) e^{-2s} z + (-2r(s) + r(2s))e^{-4s}z^2| 
\\
&\leq |r(s)| + |2r(s)| + |r(2s)| \leq C s.
\end{align*}
Since $(r_n -2 r_{n-1} +r_{n-2})/s^2$ is a difference quotient of second order, the mean value theorem gives
$$
II= s^2\sum_{n=3}^\infty r''(s_n^*)\tilde{z}^n, \qquad\mbox{for some }s_n^* \in ((n-2)s, ns).
$$
Then, by using \eqnref{rdd_gdd_equal} and the fact $|\tilde{z}|\leq e^{-s}$, we get
$$
|II| = \Bigr|s^2\sum_{n=3}^\infty \widetilde{f}''(s_n^*)\tilde{z}^n\Bigr| \leq s^2\sum_{n=3}^\infty \left|\widetilde{f}''(s_n^*)\right| e^{-ns}.
$$ 
Since $\sum_{n=3}^\infty |\widetilde{f}''(s_n^*) |e^{-ns} s$ is a Riemann sum, the RHS in the above equation can be approximated by the integral $$s\int_0^\infty |\widetilde{f}''(x)| e^{-x} dx.$$ From the assumption \eqnref{decay_fs}, we have 
\begin{align*}
|\widetilde{f}''(x)| &= 
\Big|(\p_x^2 f)(x,s) e^{2x} + 4(\p_x f)(x,s) e^{2x} + f(x,s) 4 e^{2x} \Big|
\leq Cx^N, \quad \mbox{for all }x>1
\end{align*}
and, therefore, 
$$
|II| \leq Cs \int_{0}^\infty |\widetilde{f}''(x)| e^{-x} \leq C s.
$$
So we have $|(1-\tilde{z})^2 R|=|I|+|II| \leq C s$, which implies \eqnref{def:remainderR}. The proof is completed.
\qed

Since we have $|f_0(s)|\leq C$ by \eqnref{decay_fs}, the following simplified version follows.

\begin{cor}\label{cor:simple_asymp}
Under the same assumptions on $f$ as in Theorem \ref{thm:singular_L}, the series $\mathcal{L}[f](z)$ has the following asymptotic behavior for small $s>0$:
\begin{align}
\mathcal{L}[f](z) &=  O(\frac{1}{|1-e^{-2s}z|}),
\end{align}
for $z\in \mathbb{C}$ satisfying $|z|\leq e^{s}$.
\end{cor}

\begin{remark}
We emphasize that our proposed method can be easily generalized to give any higher-order asymptotic expansions, even though only the leading order term is considered in this paper. This will be a subject of a forthcoming paper. 
\end{remark}


\subsection{Symmetric combinations of the series $\mathcal{L}$}

We shall see in section \ref{sec:stress} that the stress tensor components can be expressed in terms of the series $\mathcal{L}[f](e^{\pm\zeta+i\theta})$ for several smooth  functions $f$.
In fact, it turns out to be much more convenient to use some symmetric combinations of the series $\mathcal{L}$. The followings are those combinations. 

\begin{definition}
For a small parameter $s>0$ and a smooth function $f:\mathbb{R}^+\times (0,1) \rightarrow \mathbb{R}$, we define
\beq\label{def:M_pm}
\mathcal{M}_{\pm}[f](\zeta,\theta) := \frac{1}{2}\Big(\mathcal{L}[f](e^{\zeta + i\theta}) \pm \mathcal{L}[f](e^{-\zeta + i\theta})\Big).
\eeq
More explicitly, they can be written as
\begin{align}
\mathcal{M}_{+}[f](\zeta,\theta) &= \sum_{n=1}^\infty f(ns,s)\cosh n\zeta\, e^{in \theta}, 
\notag
\\
\mathcal{M}_{-}[f](\zeta,\theta) &= \sum_{n=1}^\infty f(ns,s)\sinh n\zeta\, e^{in \theta}. 
\label{M_pm_ch_sh}
\end{align}
\end{definition}

Then, applying Theorem \ref{thm:singular_L}, we obtain the asymptotic expansion for $\mathcal{M}_\pm[f]$ as follows.
\begin{prop}\label{lem:asymp_Msh_Mch_general_f}
We assume the same conditions on $f$ as in Theorem \ref{thm:singular_L}.
Then we have, for $|\zeta|\leq s$ and $|\theta|\leq \pi$, 
\begin{align}
\mathcal{M}_{+}[f](\zeta,\theta) &= -\frac{f_0(s)}{2}\left(1-i\frac{\sin\theta}{h(\zeta,\theta)}\right) + O(\frac{s}{h(s,\theta)}),
\label{eqn:Mpv_asymp_general_f}
\\
\mathcal{M}_{-}[f](\zeta,\theta) &= O(\frac{s}{h(s,\theta)}). 
\label{eqn:Mmv_asymp_general_f}
\end{align}
Remind that $f_0(s)=\lim_{x\rightarrow 0+}f(x,s)$.
\end{prop}
\pf
Let $\zeta\in[-s,s]$ and $\theta\in[0,\pi)$.
For notational simplicity, we introduce
\beq\label{def_z_pm}
z_\pm = e^{-2s\pm \zeta + i\theta}.
\eeq
Note that
\begin{align}
|1-z_\pm|^2&=1+e^{2(-2s\pm \zeta)} -2 e^{-2s\pm \zeta}\cos\theta 
\notag
\\
&= 2e^{-2s\pm \zeta} h(-2s\pm \zeta,\theta)
\notag
\\
&= 2e^{-2s\pm \zeta} h(\zeta,\theta) + O(s).
\label{h_estim_c}
\end{align}
It then follows 
$$
\frac{1}{|1-z_\pm|^2} \leq   \frac{C}{h(s,\theta)}.
$$
Therefore, from Theorem \ref{thm:singular_L}, we obtain
\begin{align*}
\mathcal{L}[f](e^{\pm\zeta+i\theta}) &= f_0(s) \frac{z_\pm}{1-z_\pm} + O(\frac{s}{|1-z_\pm|^2})
\\
&= f_0(s)\frac{z_\pm}{1-z_\pm} +O(\frac{s}{h(s,\theta)}).
\end{align*}
Thus, from the definition \eqnref{def:M_pm} of $\mathcal{M}_\pm$, we get
\begin{align}
\mathcal{M}_\pm [f](\zeta,\theta) 
&=  \frac{1}{2} \left(\mathcal{L}[f](e^{\zeta+i\theta}) \pm \mathcal{L}[f](e^{-\zeta+i\theta})\right ) 
\notag
\\
&= f_0(s)\frac{1}{2} \Big(  \frac{z_+}{1-z_+} \pm  \frac{z_-}{1-z_-}  \Big) +O(\frac{s}{h(s,\theta)}).
\label{M_pm_v_asymp_prev}
\end{align}
Using \eqnref{def_z_pm} and \eqnref{h_estim_c}, we compute
\begin{align*}
\frac{z_\pm }{1-z_\pm } &= \frac{z_\pm (1-\overline{z_\pm}) }{|1-z_\pm|^2}
=\frac{1}{2} \frac{\cos\theta + i\sin\theta-e^{-2s\pm \zeta}}{h( \zeta,\theta) + O(s)}
\\[.5mm]
&
=\frac{1}{2} \frac{-h(\zeta,\theta) +i \sin\theta + \cosh \zeta-e^{-2s\pm \zeta}}{h( \zeta,\theta) + O(s)} 
\\[.5mm]
&
=-\frac{1}{2} + \frac{i}{2} \frac{\sin\theta}{h(\zeta,\theta)} + O(\frac{s}{h(s,\theta)}).
\end{align*}
Therefore, the conclusion follows from \eqnref{M_pm_v_asymp_prev} and the fact that $|f_0(s)|\leq C$.
\qed


\section{Asymptotic expansion of stress tensor $\sigma$}\label{sec:stress}
In this section we derive the asymptotic expansion of the stress tensor $\sigma$ for small $s$. 
We {investigate each of the components} $(\sigma_{\theta\theta}-\sigma_{\zeta\zeta}), \sigma_{\zeta\theta}$, and $\sigma_{\zeta\zeta}$ by applying the new method of singular asymptotic expansion developed in section \ref{sec:singular}. 

\subsection{Preliminary asymptotic expansions}
{We have seen in section \ref{subsec:AiryStress_bipolar} that the stress tensor $\sigma$ is represented using the Fourier series with the coefficients $A_n$ and $B_n$ given in \eqnref{eqn:ABK}. 
To apply our new asymptotic method to those series, we need to rewrite the coefficients $A_n$ and $B_n$ using some smooth functions. 
{For doing this, we introduce some definitions. Let us define
\beq\label{def:etas_sstar}
\eta(y):=({\sinh 2y})/{(2y)}\qquad\mbox{and}\qquad
\tilde{s}: =({\sinh^2 s})/{s}.
\eeq
Note that $\eta(y) = 1+O(y^2)$ and $\tilde{s} = s+O(s^3)$.
We also define two smooth functions as follows:
\beq\label{vw}
\ds
 v(x,y):=\frac{2e^{-x}\sinh x+\eta(y) 2x}{\sinh 2x+ \eta(y) 2x}\qquad\mbox{and}\qquad w(x,y):=\frac{2x}{\sinh 2x+ \eta(y) 2x}. 
\eeq}

We have the following lemma.}
\begin{lemma}
The coefficients $A_n$ and $B_n$ in \eqnref{eqn:ABK} can be rewritten as
\begin{align}
\ds\nonumber
A_n&=\ds\frac{K}{n(n+1)}\big(v(n s,s) - \tilde{s} w(n s,s)\big), & n \geq 1,
\\[1mm] \label{newAB}
\ds B_n&=\ds \frac{K}{n(n-1)}\big(-v(n s,s) -\tilde{s} w(n s,s)\big), & n \geq 2.
\end{align}
\end{lemma}
\pf
It is easy to check that the following identity holds:
$$
ne^{\pm s}{\sinh s} = (\eta(s)  \pm \tilde{s}) ns.
$$
Then, from the definitions \eqnref{eqn:ABK} of $A_n$ and $B_n$, we have
\begin{align*}
A_n &=\ds\frac{K}{n(n+1)}
\frac{ 2(e^{-ns}\sinh ns+   (\eta(s) - \tilde{s})ns )
}{\sinh 2ns+ \eta(s)  2ns }
=\frac{K}{n(n+1)} (v(ns,s)- \tilde{s} w(ns,s)).
\end{align*}
The case of $B_n$ can be treated by the same way.
\qed

\smallskip

{
We shall see later that the stress tensor components $(\sigma_{\theta\theta}-\sigma_{\zeta\zeta}), \sigma_{\zeta\theta}$, and $\sigma_{\zeta\zeta}$ are nicely represented using $\mathcal{M}_\pm[v]$ and $\mathcal{M}_\pm[w]$. 
Therefore we need their asymptotic expansions for small $s$. We obtain the following proposition by applying Proposition \ref{lem:asymp_Msh_Mch_general_f}.
}
\begin{prop}\label{lem:asymp_Msh_Mch}
For $|\zeta|\leq s$ and $|\theta|\leq \pi$, we have
\begin{align}
\mathcal{M}_{+}[v](\zeta,\theta) &= -\frac{1}{2}\Big(1-i\frac{\sin\theta}{h(\zeta,\theta)}\Big) + O(\frac{s}{h(s,\theta)}),
\label{eqn:Mpv_asymp}
\\
\mathcal{M}_{+}[w](\zeta,\theta) &= -\frac{1}{4}\Big(1 -i\frac{\sin\theta}{h(\zeta,\theta)}\Big) + O(\frac{s}{h(s,\theta)}),
\label{eqn:Mpw_asymp}
\\
\mathcal{M}_{-}[v](\zeta,\theta) &= O(\frac{s}{h(s,\theta)}), 
\label{eqn:Mmv_asymp}
\\
  \mathcal{M}_{-}[w](\zeta,\theta) &= O(\frac{s}{h(s,\theta)}).
  \label{eqn:Mmw_asymp}
\end{align}
\end{prop}
\pf
Let us first consider the asymptotics 
\eqnref{eqn:Mpv_asymp}, \eqnref{eqn:Mmv_asymp} for $\mathcal{M}_\pm[v]$.
One can easily check that the function $v$ satisfies the conditions (A1) and (A2) in Theorem \ref{thm:singular_L}. 
We also have
$$
\lim_{x\rightarrow 0+} v(x,s) = {2}/({1+\eta(s)}) = 1 + O(s^2).
$$
Therefore, by applying Proposition \ref{lem:asymp_Msh_Mch_general_f} to $\mathcal{M}_{\pm}[v]$, we immediately get \eqnref{eqn:Mpv_asymp} and \eqnref{eqn:Mmv_asymp}.

The asymptotics \eqnref{eqn:Mpw_asymp}, \eqnref{eqn:Mmw_asymp} for $\mathcal{M}_\pm[w]$ can be proved by the exactly same way using the fact that 
$$
\lim_{x\rightarrow 0+} w(x,s) = {1}/{(1+\eta(s))} = 1/2 + O(s^2).
$$
The proof is completed.
\qed

Now we are ready to derive the asymptotic expansion of the stress tensor.

\subsection{Asymptotic expansion of $\sigma_{\theta\theta}-\sigma_{\zeta\zeta}$}\label{sec:asymp1}

In this subsection, we represent $\sigma_{\theta\theta}^1-\sigma_{\zeta\zeta}^1$ in terms of $\mathcal{M}_\pm$ and then derive its asymptotic expansion for small $s$.

We have the following lemma whose proof will be given in Appendix \ref{sec:stress_Mpm}.

\begin{lemma} \label{lem:sigma_tt_zz}
The stress component $\sigma_{\theta\theta}^1-\sigma_{\zeta\zeta}^1$ can be represented using $\mathcal{M}_\pm[v]$ and $\mathcal{M}_\pm[w]$ as follows: for $|\zeta|\leq s$ and $|\theta|\leq \pi$,
\begin{align}
\ds(\sigma_{\theta\theta}^1-\sigma_{\zeta\zeta}^1)(\zeta,\theta)
&=
4Kh(\zeta,\theta)  \sinh\zeta\,\mathrm{Re}\big\{ \mathcal{M}_{-}[v](\zeta,\theta)\big\}
\notag
\\
&\quad
-4Kh(\zeta,\theta) \tilde{s}\cosh\zeta\, \mathrm{Re}\big\{\mathcal{M}_{+}[w](\zeta,\theta) \big\}
\notag
\\
&\quad
+2K \sinh^2 \zeta.
\label{eqn:sigma_tt_zz}
\end{align}

\end{lemma}

We have the following asymptotic result for $\sigma^1_{\theta\theta}-\sigma^1_{\zeta\zeta}$.

\begin{prop} \label{prop:asymp_sigma_ttzz}
For small $s$, we have the asymptotic expansion of $\sigma^1_{\theta\theta}-\sigma^1_{\zeta\zeta}$ as follows:
for $|\zeta|\leq s$ and $|\theta|\leq \pi$,
$$
(\sigma^1_{\theta\theta}-\sigma^1_{\zeta\zeta})(\zeta,\theta) =  
\frac{1}{ \mathcal{I}_0 s} h(\zeta,\theta)+O(1).
$$

\end{prop}

\pf 
We shall apply Proposition \ref{lem:asymp_Msh_Mch}
to \eqnref{eqn:sigma_tt_zz}. 
Let us begin with the first term in RHS of \eqnref{eqn:sigma_tt_zz}. Remind that $K=O(s^{-2})$ and $\zeta = O(s)$. 
So, by applying Proposition \ref{lem:asymp_Msh_Mch} and \eqnref{s2_h_bdd}, we obtain
\begin{align*}
\Big|4Kh(\zeta,\theta)  \sinh\zeta\,\mbox{Re}\big\{ \mathcal{M}_{-}[v](\zeta,\theta)\big\}\Big|
 \leq C s^{-2} h(s,\theta) \,s\frac{s}{h(s,\theta)}  \leq C.
\end{align*}
We next consider the second term. By Lemma \ref{lem:asymp_Ks} and Proposition \ref{lem:asymp_Msh_Mch}, we have
\begin{align*}
4Kh(\zeta,\theta)  & \tilde{s}\cosh\zeta\,\mbox{Re}\big\{ \mathcal{M}_{+}[w](\zeta,\theta)\big\}
\\
&=
4 K \tilde{s} \cosh\zeta  h(\zeta,\theta) (-\frac{1}{4}+ O(\frac{s}{h(s,\theta)}))
\\
&=
 -K \tilde{s}\cosh\zeta h(\zeta,\theta) + O(1)
 \\
 &
 =(-\frac{1}{s^2 \mathcal{I}_0} + O(s^{-1})) (s+O(s^3)) h(\zeta,\theta)  + O(1)
\\
&=-\frac{1}{ \mathcal{I}_0 s} h(\zeta,\theta)+O(1).
\end{align*}
For the third term, it is clear from Lemma \ref{lem:asymp_Ks} that
$2K\sinh^2\zeta = O(1)$. The proof is completed.
\qed

\subsection{Asymptotic expansion of $\sigma_{\zeta\theta}$}
We consider the stress component $\sigma_{\zeta\theta}$.
We have the following lemma whose proof will be given in Appendix \ref{sec:stress_Mpm}.

\begin{lemma}\label{lem:sigma_zt}
The stress component $\sigma_{\zeta\theta}^1$ can be represented using $\mathcal{M}_{\pm}[v]$ and $\mathcal{M}_{\pm}[w]$ as follows: for $|\zeta|\leq s$ and $|\theta|\leq \pi$,
\begin{align}
\ds \sigma_{\zeta\theta}^1(\zeta,\theta)
&=
2Kh(\zeta,\theta)  \sinh\zeta\,\mathrm{Im}\big\{ \mathcal{M}_{+}[v](\zeta,\theta)\big\}
\notag
\\
&\quad
-2Kh(\zeta,\theta) \tilde{s}\cosh\zeta\, \mathrm{Im}\big\{\mathcal{M}_{-}[w](\zeta,\theta) \big\}
\notag
\\
&\quad
-K \sinh \zeta \sin\theta.
\label{eqn:sigma_zt}
\end{align}

\end{lemma}

We have the following asymptotic result for $\sigma_{\zeta\theta}$.

\begin{prop}\label{prop:asymp_sigma_zt}
For small $s>0$, the stress component $\sigma_{\zeta\theta}^1$ has following asymptotic behavior:
for $|\zeta|\leq s$ and $|\theta|\leq \pi$,
$$
\sigma_{\zeta\theta}^1(\zeta,\theta) =  O(1).
$$

\end{prop}
\pf
We apply Proposition \ref{lem:asymp_Msh_Mch} to \eqnref{eqn:sigma_zt} as in the proof of Proposition \ref{prop:asymp_sigma_ttzz}. 

Let us consider the first term in RHS of \eqnref{eqn:sigma_zt}. 
Remind that $K=O(s^{-2}), \tilde{s}=O(s)$ and $\zeta=O(s)$.
By Proposition \ref{lem:asymp_Msh_Mch} and \eqnref{s2_h_bdd}, we have
\begin{align}
\notag 2Kh(\zeta,\theta) & \sinh\zeta\,\mathrm{Im}\big\{ \mathcal{M}_{+}[v](\zeta,\theta)\big\} 
\\
\notag &= 2K h(\zeta,\theta) \sinh\zeta \Big(\frac{1}{2}\frac{\sin\theta}{h(\zeta,\theta)} + O(\frac{s^2}{[h(s,\theta)]^{3/2}})\Big)
\\
&=K\sinh\zeta {\sin\theta} +O(1).
\label{pf_asymp_sigma_zt1}
\end{align}

Next we consider the second term. 
By Proposition \ref{lem:asymp_Msh_Mch} and \eqnref{s2_h_bdd}, we have
\beq\label{pf_asymp_sigma_zt2}
\Big|2Kh(\zeta,\theta) \tilde{s}\cosh\zeta\, \mathrm{Im}\big\{\mathcal{M}_{-}[w](\zeta,\theta) \big\}\Big|
\leq C s^{-2} h(s,\theta) s \frac{s}{h(s,\theta)} \leq C.
\eeq
From \eqnref{eqn:sigma_zt}, \eqnref{pf_asymp_sigma_zt1} and \eqnref{pf_asymp_sigma_zt2}, the conclusion immediately follows.
\qed

\subsection{Asymptotic expansion of $\sigma_{\zeta\zeta}$}
We consider the stress component $\sigma_{\zeta\zeta}$. 
We will represent the stress component $\sigma_{\zeta\zeta}^1$ using $\widetilde{\mathcal{M}}_{\pm}$ (instead of ${\mathcal{M}}_{\pm}$) defined by
\beq\label{def:widetilde_M_pm}
\widetilde{\mathcal{M}}_{\pm}[f](\zeta,\theta) := \frac{1}{2}{\rm{Re}}\left\{\mathcal{L}[f](e^{\zeta + i\theta})e^{2\zeta + 2i\theta} \pm \mathcal{L}[f](e^{-\zeta + i\theta})e^{-2\zeta + 2i\theta}\right\}.
\eeq
One can easily check that
\begin{align}
\widetilde{\mathcal{M}}_{+}[f](\zeta,\theta) = \sum_{n=3}^\infty f\big((n-2)s,s\big)\cosh n\zeta\, \cos{n\theta},
\notag 
\\
\widetilde{\mathcal{M}}_{-}[f](\zeta,\theta) = \sum_{n=3}^\infty f\big((n-2)s,s\big)\sinh n\zeta\, \cos{n \theta}.
\label{widetilde_M_pm_ch_sh}
\end{align}
We have the following lemma for $\sigma_{\zeta\zeta}^1$. See Appendix \ref{sec:stress_Mpm} for its proof.

\begin{lemma}\label{lem:sigma_zz_series}
For $|\zeta|\leq s$ and $|\theta|\leq \pi$, the stress component $\sigma_{\zeta\zeta}^1$ has the following representation:
\begin{align}
\sigma_{\zeta\zeta}^1(\zeta,\theta) &=  
 2K s \sinh^2\zeta \,\widetilde{\mathcal{M}}_+ [v_1]
+ K s^2 \tilde{s} \cosh^2\zeta \,\widetilde{\mathcal{M}}_+[w_2]
\notag
\\
&\quad +  (K s^2/2)\sinh 2\zeta \,\widetilde{\mathcal{M}}_-[v_2]
+ K s \tilde{s} \sinh 2\zeta \,\widetilde{\mathcal{M}}_-[w_1]
\notag
\\
&\quad - K s^3 \widetilde{\mathcal{M}}_+[v_3] 
+ 2K s^2 \tilde{s} \widetilde{\mathcal{M}}_+[w_3] 
+O(1)
\label{sigma_zz_tilde_M_pm}
\end{align}
where the functions $v_j$ and $w_j$ are given by
\begin{align*}
&v_1(x):=\frac{v(x+3s,s)-v(x+s,s)}{2s}
,\quad
 v_2(x) :=
\frac{v(x+3s,s)-2v(x+2s,s) + v(x+s,s)}{s^2},
\\
& w_1(x):=\frac{w(x+3s,s)-w(x+s,s)}{2s},
\quad 
w_2(x) :=
\frac{w(x+3s,s)-2w(x+2s,s) + w(x+s,s)}{s^2},
\\
& v_3(x) :=v_2(x)/(x+2s),
\quad
w_3(x) :=w_1(x)/(x+2s).
\end{align*}
Note that $v_1,v_2,w_1$ and $w_2$ are difference quotients of $v$ or $w$.

\end{lemma}

We have the following asymptotic result for $\sigma_{\zeta\zeta}^1$.
\begin{prop}\label{prop:asymp_sigma_zz}
The stress component $\sigma_{\zeta\zeta}^1$ has the following asymptotic behavior for small $s>0$:
for $|\zeta|\leq s$ and $|\theta|\leq \pi$,
$$
\sigma_{\zeta\zeta}^1 (\zeta,\theta)=  
O(1).
$$
\end{prop}
\pf
Using the mean value thoerem, it is easy to check that $v_1,v_2,w_1$ and $w_2$ satisfy the condition \eqnref{decay_fs}.  Moreover, $v_2$ and $w_1$ satisfy 
$$|(v_2)^{(k)}(x)| + |(w_1)^{(k)}(x)| \leq C x, \quad \mbox{for small } x>0 \mbox{ and } k = 1,2,3.  
$$
Using the above estimates, it is also easy to check that $v_3$ and $w_3$ satisfy \eqnref{decay_fs}.
Therefore, we can apply our asymptotic method to $\widetilde{\mathcal{M}}_\pm[v_j]$ and $\widetilde{\mathcal{M}}_\pm[w_j]$ for all $j=1,2,3$.
By using Corollary \ref{cor:simple_asymp} and \eqnref{def:widetilde_M_pm}, we have
$$
\widetilde{\mathcal{M}}_\pm[v_j], \widetilde{\mathcal{M}}_\pm[w_j] = O({ |1-e^{-s+i\theta}|^{-1} })= O(s^{-1}), \qquad j=1,2,3.
$$
Then we can show that each of all terms in RHS of \eqnref{sigma_zz_tilde_M_pm} is of $O(1)$. For simplicity, we consider the first term only. Remind that $K=O(s^{-2})$ and $\zeta = O(s)$. We have
$$
|2K s \sinh^2 \zeta \, \mathcal{M}_+[v_1]| \leq C s^{-2} s s^2 s^{-1} \leq C.
$$
The other terms can be estimated in the exactly same way. The proof is completed.
\qed

\subsection{Asymptotic expansion of the stress tensor $\sigma$}\label{subsec:asymp_stress_tensor}

In this subsection, we finally derive an asymptotic expansion of the stress tensor $\sigma$ for small $s$.

As an immediate consequence of Propositions \ref{prop:asymp_sigma_ttzz}, \ref{prop:asymp_sigma_zt}, and \ref{prop:asymp_sigma_zz}, we get the following proposition.
\begin{prop} 
Let $\mathbf{u}$ be the solution to the two circular holes problem \eqnref{elas_eqn}. 
For small $s>0$,
we have the following asymptotic expansion of the corresponding stress tensor $\sigma$:
for $|\zeta|\leq s$ and $|\theta|\leq \pi$,
\beq\label{stress_tensor_final_asymp_bipolar}
\sigma(\zeta,\theta) = \frac{1}{\mathcal{I}_0 s}h(\zeta,\theta) \mathbf{e}_{\theta} \otimes \mathbf{e}_{\theta} + O(1),
\eeq
where $h(\zeta,\theta)$ is defined as \eqnref{scale} and $\mathcal{I}_0$ is given by \eqnref{def:I_0}.
\end{prop}

Now we are ready to prove Theorem \ref{thm:main_stress_asymp}, which is the main result in this paper.

\smallskip

\noindent{\sl Proof of Theorem \ref{thm:main_stress_asymp}}.
It is sufficient to rewrite \eqnref{stress_tensor_final_asymp_bipolar} in a coordinate-free form. 
From \eqnref{unit_basis} and \eqnref{grad_f_bipolar}, we have
$$
h(\zeta,\theta) = \alpha|\nabla \zeta|, \quad \mathbf{e}_\theta = \nabla\theta/|\nabla \theta|.
$$
Therefore, since $\alpha=\sqrt{r\ep}+O(\ep^{3/2})$ and $s=\sqrt{\ep/ r} + O(\ep^{3/2})$ as $\ep\rightarrow 0$, we obtain
\begin{align*}
\sigma &= \frac{1}{\mathcal{I}_0 }\frac{\alpha}{s} \frac{h(\zeta,\theta)}{\alpha} \mathbf{e}_{\theta} \otimes \mathbf{e}_{\theta} +O(1)
\\
&= \frac{r}{\mathcal{I}_0 } |\nabla \zeta| (\frac{\nabla \theta}{|\nabla \theta|}\otimes \frac{\nabla \theta}{|\nabla \theta|}) + O(1).
\end{align*}
 
Now we represent $\nabla \zeta$ and $\nabla \theta$ in a coordinate free form.
Using $\mathbf{p}_1 = (-\alpha,0), \mathbf{p}_2 = (\alpha,0)$ and \eqnref{xitheta1}, we can easily see that
\begin{align*}
\nabla \zeta &= \frac{(x+\alpha,y)}{(x+\alpha)^2+y^2} - \frac{(x-\alpha,y)}{(x-\alpha)^2+y^2} 
=\frac{\Bx-\Bp_1}{|\Bx-\Bp_1|^2} -\frac{\Bx-\Bp_2}{|\Bx-\Bp_2|^2},
\\[0.3em]
\nabla \theta &= \frac{(-y,x+\alpha)}{(x+\alpha)^2+y^2} - \frac{(-y,x-\alpha)}{(x-\alpha)^2+y^2} 
=\frac{(\Bx-\Bp_1)^\perp}{|\Bx-\Bp_1|^2} - \frac{(\Bx-\Bp_2)^\perp}{|\Bx-\Bp_2|^2}.
\end{align*}
The proof is completed.
\qed




\appendix

\section{Proof of Lemma \ref{lem:asymp_Ks}}\label{appendix:Ks}
Here we derive the asymptotic expansion of the constant $K$ for small $s>0$. 
In \cite{LM}, the asymptotics of the constant $K$ was already derived. 
However, we give a simpler method than theirs. 
We apply the following summation formula. 
\begin{lemma}(Euler-Maclaurin formula)\label{EM} Let $N\in\mathbb{N}$ and let $f \in C^N(\mathbb{R}^+)\cap L^1(\mathbb{R}^+)$. Then, for a small parameter $s>0$,  
we have
\begin{align}
\sum_{n=0}^\infty f(a+n s) s = \int_a^\infty f(x) dx +\frac{s}{2}f(a)-
\sum_{m=2}^N s^m \frac{B_m}{m!} f^{(m-1)}(a) + R_N(a)
\end{align}
where $B_m$ is the Bernoulli numbers and the remainder term $R_N$ satisfies
\beq
\left|R_N(a)\right| \leq 4 \left( \frac{s}{2\pi}\right)^N \int_a^\infty \left|f^{(N)}(x)\right|dx.
\eeq
\end{lemma}

\noindent\textbf{Proof of Lemma \ref{lem:asymp_Ks}.}
Recall that
\beq\label{def_K_App}
K= \big(\frac{1}{2}+\tanh{s}\sinh^2{s}-  4 P(s) \big)^{-1},
\eeq
where $P(s)$ is given by
\beq\label{def_P_App}
P(s)=
\sum_{n=2}^{\infty}\frac{e^{-n s}\sinh ns+n(n\sinh s+\cosh s)\sinh s}{n(n^2-1)(\sinh 2ns + n\sinh 2s )}.
\eeq

One can easily check that 
$$
\lim_{s\rightarrow 0}P(s)=\sum_{n=2}^\infty \frac{1}{2n(n^2-1)}=\frac{1}{8}.
$$
In view of this, we decompose $P(s)$ as
\begin{align}
P(s)&=
\frac{1}{8} + \bigg(P(s)  - \sum_{n=2}^\infty \frac{1}{2n(n^2-1)}\bigg)
\nonumber 
\\
&=\frac{1}{8}-s^2\sum_{n=2}^\infty \frac{
\sinh^2 ns - n^2\sinh^2 s}{ ns\big((ns)^2-s^2\big)(\sinh 2ns+ n \sinh 2s)}s\nonumber\\
&=\frac{1}{8} - s^2 \sum_{n=2}^\infty f_K(ns)s,\label{P_re}
\end{align}
where
$$
f_{K}(x)= \frac{\sinh^2 x-(\sinh^2s/s^2)x^2}{x(x^2-s^2)(\sinh 2 x+ 2\eta(s) x)}, \quad x>s.
$$

One can easily see that
$$
|f_K(2s)| \leq C.
$$
Straightforward but tedious computations give us
$$
  |f_K'(x)|  \leq C {(1+|x|)^{-3}}, \quad\mbox{for } 2s\leq x <\infty.
$$
So we have
$$
\int_{2s}^\infty |f_K'(x)| dx \leq C.
$$
Remind that $C$ denotes a positive constant independent of $s$.

Therefore, by applying Lemma \ref{EM} to $f_K$ with $a=2s$, we obtain
\begin{align*}
\sum_{n=2}^\infty f_K(ns)s 
&=\int_{2s}^\infty f_K(x) dx +\frac{s}{2}f_K(2s) + O\Big(s \int_{2s}^\infty |f_K'(x)| dx \Big)
\\
&= \int_{2s}^\infty f_K(x) dx +O(s).
\end{align*}

We now make an approximation of $f_K$ for small $s$. One can see that $f_K$ satisfies
$$
\big|f_K(x)-f_{0}(s)\big| \leq C s^2 (1+|x|)^{-5}\quad\mbox{for }\quad x\geq 2s,$$
where $f_0$ is defined by $$f_0(x)=\frac{\sinh^2 x-x^2}{x^3 (\sinh 2x+ 2x)}, \quad x>0.$$
Note that $f_0$ has a removable singularity at $x=0$.

Therefore, we obtain
$$
\sum_{n=2}^\infty f_K(ns)s = \int_{2s}^\infty f_0(x) dx +O(s) = \int_{0}^\infty f_0(x) dx +O(s).
$$
So, from \eqnref{P_re}, we get
\beq\label{P_asymp_app}
P(s) = \frac{1}{8} - s^2 \int_0^\infty f_0(x) dx + O(s^3).
\eeq

Now we are ready to get the asymptotics of $K$. From \eqnref{def_K_App} and \eqnref{P_asymp_app}, we obtain
\begin{align*}
K&=\big(\frac{1}{2}+O(s^3)-  4 P(s) \big)^{-1} 
= \Big( 4 s^2\int_0^\infty f_0(x) dx + O(s^3) \Big)^{-1} 
\\
&= \frac{1}{s^2} \frac{1}{4}\Big( \int_0^\infty f_0(x) dx\Big)^{-1} + O(s^{-1}).
\end{align*}
The proof is completed.
\qed

\section{Proofs of Lemmas \ref{lem:sigma_tt_zz}, \ref{lem:sigma_zt} and \ref{lem:sigma_zz_series} }\label{sec:stress_Mpm}

\noindent\textbf{Proof of Lemma \ref{lem:sigma_tt_zz}.}
We obtain, from \eqnref{eqn:sigma1} and \eqnref{newAB}, that
\begin{align}
(\sigma^1_{\theta\theta}-\sigma^1_{\zeta\zeta})(\zeta,\theta)
\notag&= 
K\big(\cosh 2\zeta-2\cosh\zeta \cos\theta+\cos 2\theta\big) + 2K h(\zeta,\theta) \big(v(s,s)+\tilde{s}w(s,s)\big) \cos\theta
\\
\notag&\quad+  2K h(\zeta,\theta)\sum_{n=1}^{\infty} 
v(ns,s)\big[\cosh(n+1)\zeta - \cosh(n-1)\zeta\big] \cos n\theta 
\\
\notag&\quad -2K h(\zeta,\theta) \tilde{s}\sum_{n=1}^{\infty} 
w(ns,s)\big[\cosh(n+1)\zeta + \cosh(n-1)\zeta\big] \cos n\theta. 
\end{align}
{One can easily check that $v(s,s)+\tilde{s}w(s,s)=1$. Recall the hyperbolic identities
\begin{align*}
\cosh(n+1)\zeta + \cosh(n-1)\zeta &= 2\cosh \zeta\cosh n\zeta, 
\\
\cosh(n+1)\zeta - \cosh(n-1)\zeta &= 2\sinh \zeta\sinh n\zeta.
\end{align*}
Then, in view of the expressions \eqnref{M_pm_ch_sh} for $\mathcal{M}_\pm$, we see that
\begin{align}
(\sigma^1_{\theta\theta}-\sigma^1_{\zeta\zeta})(\zeta,\theta)
\notag&= 
K\big(\cosh 2\zeta+\cos 2\theta - 2\cos^2\theta \big)
\\
&\quad 
+2Kh(\zeta,\theta)  \sinh\zeta\,\mathrm{Re}\big\{ \mathcal{M}_{-}[v](\zeta,\theta)\big\}
\notag
\\
&\quad
-2Kh(\zeta,\theta) \tilde{s}\cosh\zeta\, \mathrm{Re}\big\{\mathcal{M}_{+}[w](\zeta,\theta) \big\}.
\notag
\end{align}
Since 
$$\cosh 2\zeta + \cos 2\theta - 2\cos^2\theta = \cosh2\zeta -1 = \sinh^2\zeta,$$ 
we get the conclusion.}
\qed

\smallskip

\noindent\textbf{Proof of Lemma \ref{lem:sigma_zt}.}
We obtain, from \eqnref{eqn:sigma2} and \eqnref{newAB}, we have
\begin{align*}
\sigma_{\zeta\theta}^1(\zeta,\theta)&=-K\sinh\zeta\sin\theta 
+  K h(\zeta,\theta)\sum_{n=1}^{\infty} 
v(ns,s)\big[\sinh(n+1)\zeta - \sinh(n-1)\zeta\big] \sin n\theta 
\\
\notag&\quad -K h(\zeta,\theta) \tilde{s}\sum_{n=1}^{\infty} 
w(ns,s)\big[\sinh(n+1)\zeta + \sinh(n-1)\zeta\big] \sin n\theta. 
\end{align*}
{Then, from the hyperbolic identities 
\begin{align*}
\sinh(n+1)\zeta + \sinh(n-1)\zeta &= 2\cosh \zeta\sinh n\zeta, 
\\
\sinh(n+1)\zeta - \sinh(n-1)\zeta &= 2\sinh \zeta\cosh n\zeta,
\end{align*} }
and the expressions \eqnref{M_pm_ch_sh} for $\mathcal{M}_\pm$, the conclusion follows.
\qed

\smallskip

\noindent\textbf{Proof of Lemma \ref{lem:sigma_zz_series}.}
{From \eqnref{eqn:sigma3} and the trigonometric identities 
\begin{align*}
2\cos \theta \cos n\theta &= \cos(n+1)\theta + \cos(n-1)\theta,
\\
2\sin \theta \sin n\theta &= -\cos(n+1)\theta + \cos(n-1)\theta,
\end{align*} }
we have
\begin{align*}
\ds\sigma_{\zeta\zeta}^1(\zeta,\theta)
&\ds=-\frac{K}{2}(\cosh2\zeta-2\cosh \zeta \cos\theta +\cos 2\theta)
\ds+\frac{1}{2}\sum_{n=1}^\infty\Big((n^2+n)\phi_n(\zeta)\cos(n-1)\theta \\
\ds&\quad  -2\Big[(n^2-1)\cosh\zeta\,\phi_n(\zeta)+\sinh\zeta\,\phi_n'(\zeta)\Big]\cos n\theta + (n^2-n)\phi_n(\zeta)\cos(n+1)\theta\Big).
\end{align*}
Then, by translating summation indices, we get 
\beq\label{sigma_zz_p_sum_psi_n}
\sigma_{\zeta\zeta}^1(\zeta,\theta)=p(\zeta,\theta)+\frac{1}{2}\sum_{n=3}^\infty\psi_n(\zeta)\cos n\theta,
\eeq
where $\psi_n(\zeta)$ and $p(\zeta,\theta)$ are given by 
\begin{align}
\ds \psi_n(\zeta)&=
(n+1)(n+2)\phi_{n+1}(\zeta) - 2(n^2-1)\cosh\zeta\,\phi_n(\zeta)
\notag
\\
&\qquad
+(n-1)(n-2)\phi_{n-1}(\zeta)-2\sinh\zeta\,\phi_n'(\zeta), \qquad n\geq 1,
\label{def:psi_n}
\\
\ds
p(\zeta,\theta)&=-(K/2)(\cosh2\zeta-2\cosh \zeta \cos\theta +\cos 2\theta)
\notag
\\&\qquad 
+ \phi_1(\zeta)+ (1/2)[\psi_1(\zeta)\cos\theta + \psi_2(\zeta)\cos2\theta],
\label{def:p_zz}
\end{align}
with $\phi_0(\zeta) =0$. 

Let us estimate $p(\zeta,\theta)$. Remind that $K=O(s^{-2})$ and $\zeta = O(s)$. 
One can easily see that $\phi_1(\zeta)=K/2 + O(1)$, $\psi_1(\zeta) = -2K+O(1)$ and $\psi_2(\zeta) = K +O(1)$. So we have
\begin{align}
p(\zeta,\theta) &= -(K/2)\cosh 2\zeta + \phi_1(\zeta) + \Big[ K\cosh\zeta +\psi_1(\zeta)/2 \Big]\cos\theta
\notag
\\
&\quad +
\Big[ -K/2 +\psi_2(\zeta)/2 \Big]\cos 2\theta
=O(1).
\label{p_estim}
\end{align}

We now compute $\psi_n$ in terms of $v$ and $w$.  
From \eqnref{phin}, \eqnref{newAB} and the hyperbolic identities, we have
\begin{align*}
(n+1)(n+2)\phi_{n+ 1}(\zeta) &= 
(n+1)(n+2) A_{n+1} \cosh(n+2)\zeta 
+
(n+1)(n+2) B_{n+1} \cosh n\zeta
\\
&=
K v_{n+1} (\cosh(n+2)\zeta - \cosh n\zeta)
-
K \tilde{s}w_{n+1} (\cosh(n+2)\zeta + \cosh n\zeta)
\\&
\quad-
\frac{2 K}{n} (v_{n+1} +\tilde{s} w_{n+1})\cosh n \zeta,
\\
&=
2 Kv_{n+1} \sinh^2 \zeta \cosh n\zeta
-2 K \tilde{s}  w_{n+1} \cosh^2\zeta \cosh n\zeta
\\
&\quad +K(v_{n+1} - \tilde{s}w_{n+1}) \sinh 2\zeta \sinh n\zeta
-\frac{2K}{n} (v_{n+1} +\tilde{s} w_{n+1})\cosh n \zeta,
\end{align*}
where $v_n=v(ns,s)$ and $w_n=w(ns,s)$.
Similarly, we have
\begin{align*}
(n-1)(n-2)\phi_{n- 1}(\zeta) 
&=
-2K v_{n-1}\sinh^2 \zeta \cosh n\zeta
-2K \tilde{s} w_{n-1}\cosh^2 \zeta \cosh n\zeta%
\\
&\quad +K(v_{n-1} + \tilde{s}w_{n-1}) \sinh 2\zeta \sinh n\zeta
-\frac{2K}{n} (v_{n+1} -\tilde{s} w_{n+1})\cosh n \zeta
,
\\
(n^2-1)\phi_n(\zeta) &= {2K \big(v_n + \frac{\tilde{s} w_n}{n} \big)}\sinh \zeta \sinh n \zeta -{2K (\tilde{s}w_n + \frac{v_n}{n}) } \cosh\zeta \cosh n\zeta,
\\
\phi_n'(\zeta) &= \frac{2K}{n} v_n \sinh \zeta \cosh n \zeta - \frac{2K}{n} \tilde{s} w_n \cosh \zeta \sinh n \zeta.
\end{align*}
By substituting these expressions into \eqnref{def:psi_n}, we get 
\begin{align*}
\frac{1}{K}\psi_n(\zeta) &= 
\big[v_{n+1} -v_{n-1} \big]2\sinh^2\zeta \cosh n\zeta +
 \big[w_{n+1}-2 w_n +w_{n-1} \big] 2\tilde{s}\cosh^2\zeta \cosh n\zeta
 \notag
\\
&\quad + \big[  (v_{n+1}-2v_n + v_{n-1})  +  \tilde{s}(w_{n+1}-w_{n-1})\big] \sinh 2\zeta \sinh n\zeta
\notag
\\
&\quad +\big[-
(v_{n+1}-2 v_n +v_{n-1} )/n +  \tilde{s} (w_{n+1}-w_{n-1})/n \big]2\cosh n\zeta.
\end{align*}
Then, from \eqnref{widetilde_M_pm_ch_sh}, \eqnref{sigma_zz_p_sum_psi_n} and \eqnref{p_estim}, the conclusion follows.
\qed

\section*{Acknowledgements}

We would like to thank Graeme W. Milton for pointing out to us the existence of reference \cite{MM}.


\begin{thebibliography}{99}

\bibitem{ACKLY-ARMA-13} H. Ammari, G. Ciraolo, H. Kang, H. Lee and K. Yun, Spectral analysis of the Neumann-Poincar\'e operator and characterization of the stress concentration in anti-plane elasticity, Arch. Ration. Mech. Anal. 208 (2013), 275--304.


\bibitem{AKLLL-JMPA-07} H. Ammari, H. Kang, H. Lee, J. Lee and M. Lim, Optimal bounds on the gradient of solutions to conductivity problems, J. Math. Pure. Appl. 88 (2007), 307--324.

\bibitem{AKLLZ-JDE-09} H. Ammari, H. Kang, H. Lee, M. Lim and H. Zribi, Decomposition theorems and fine estimates for electrical fields in the presence of closely located circular inclusions, J. Differ. Equations 247 (2009), 2897-2912.

\bibitem{AKL-MA-05} H. Ammari, H. Kang and M. Lim, Gradient estimates for solutions to the conductivity problem, Math. Ann. 332(2) (2005), 277--286.



\bibitem{BLY-ARMA-09}  E.S. Bao, Y. Li, and B. Yin,
Gradient estimates for the perfect conductivity problem, Arch. Rat. Mech. Anal. 193 (2009), 195-226.

\bibitem{BLL-ARMA-15} J. Bao, H. Li and Y. Li, Gradient Estimates for Solutions of the Lam\'e
System with Partially Infinite Coefficients, Arch. Ration. Mech. Anal. 215 (2015) 307--351.

\bibitem{BLL-AM-17} J. Bao, Y. Li and H. Li,
Gradient estimates for solutions of the Lame system with partially infinite coefficients in dimensions greater than two, Adv. Math. 305 (2017), 298-338.

\bibitem{BLY-CPDE-10} E.S. Bao, Y. Li, and B. Yin,
Gradient estimates for the perfect and insulated conductivity problems with multiple inclusions, Commun. Part. Diff. Eq. 35 (2010), 1982--2006.


\bibitem{CM88}  C. Callias and X. Markenscoff, Singular asymptotics of integrals and the near-field radiated from nonuniformly moving dislocations, Arch. Ration. Mech. Anal. 102 (1988), 273--285.

\bibitem{CM93} C. Callias and X. Markenscoff, The singularity of the stress field of two nearby holes in a planar elastic medium, Quarterly of Applied Mathematics 51 (1993), 547--557.

\bibitem{Gorb-MMS-16} Y. Gorb, Singular behavior of electric field of high contrast concentrated composites, SIAM Multi. Model. Simul. 13(4) (2015), 1312–-1326.

\bibitem{GN-MMS-12} Y. Gorb and A. Novikov, Blow-up of solutions to a p-Laplace equation, SIAM Multi. Model. Simul. 10 (2012), 727--743.



\bibitem{Jeff} G.B. Jeffery,
Plane Stress and Plane Strain in Bipolar Coordinates, Philosophical Transactions of the Royal Society of London. Series A, Containing Papers of a Mathematical or Physical Character, Vol. 221 (1921), 265-293.

\bibitem{KLY-MA-15} H. Kang, H. Lee and K. Yun, Optimal estimates and asymptotics for the stress concentration between closely located stiff inclusions, Math. Annalen 363 (2015), 1281--1306.


\bibitem{KK-CM-16}
H. Kang and E. Kim, Estimation of stress in the presence of closely located elastic inclusions: A numerical study, Contemporary Math. 660 (2016), 45--57.


\bibitem{KLY-JMPA-13} H. Kang, M. Lim and K. Yun, Asymptotics and computation of the solution to the conductivity
equation in the presence of adjacent inclusions with extreme conductivities, {J. Math. Pure. Appl.} 99 (2013), 234--249.

\bibitem{KLY-SIAP-14} H. Kang, M. Lim and K. Yun, Characterization of the electric field concentration between two adjacent spherical perfect conductors, SIAM J. Appl. Math. 74 (2014), 125--146.




\bibitem{Lekner-PRSA-12} J. Lekner, Electrostatics of two charged conducting spheres, {Proc. R. Soc. A,} 468 (2012), 2829--2848.


\bibitem{CBLing} C. B. Ling, On the Stresses in a Plate Containing Two Circular Holes, Journal of Applied Physics, 19, 77 (1948)

\bibitem{LLBY-QAM-14} H. Li, Y. Li, E.S. Bao and B. Yin,
Derivative estimates of solutions of elliptic systems in narrow domains, Quarterly of Applied Mathematics 72 (2014), 589--596.

\bibitem{LY-JMAA-15} M. Lim and S. Yu, Asymptotics of the solution to the conductivity equation in the presence of adjacent circular inclusions with finite conductivities, J. Math. Anal. Appl. 421 (2015), 131--156.

\bibitem{LY-CPDE-09} M. Lim and K. Yun, Blow-up of electric fields between closely spaced spherical perfect conductors,
Commun. Part. Diff. Eq. 34 (2009), 1287--1315.

\bibitem{LY-JDE-11} M. Lim and K. Yun, Strong influence of a small fiber on shear stress in fiber-reinforced composites,
J. Differ. Equations 250 (2011), 2402--2439.

\bibitem{LN-CPAM-03}
Y.Y. Li and L. Nirenberg, Estimates for elliptic systems from composite material. Comm. Pure Appl. Math. 56 (2003), 892--925.

\bibitem{LV-ARMA-00}
Y.Y. Li and M. Vogelius, Gradient estimates for solutions to divergence form elliptic equations with discontinuous coefficients. Arch. Ration. Mech. Anal. 135 (2000), 91--151.



\bibitem{MM} R. C. McPhedran and A. B. Movchan, The Rayleigh multipole method for linear elasticity, J. Mech. Phys. Solids 42(5) (1994), 711-727.


\bibitem{MPM} R. C. McPhedran, L. Poladian, and G. W. Milton, Asymptotic studies of closely spaced, highly conducting cylinders, Proceedings of the Royal Society of London A: Mathematical, Physical and Engineering Sciences, Vol. 415. No. 1848 (1988), 185-196.
\bibitem{Markenscoff} X. Markenscoff, Stress amplification in vanishingly small geometries, Computational Mechanics, 19(1) (1996), 77-83.




\bibitem{WM96} L. Wu and X. Markenscoff, Singular stress amplification between two holes in tension, Journal of Elasticity, 44(2) (1996), 131-144.


\bibitem{Yun-SIAP-07} K. Yun, Estimates for electric fields blown up between closely adjacent conductors with arbitrary shape, SIAM J. Appl. Math. 67 (2007), 714--730.

\bibitem{Yun-JMAA-09} K. Yun, Optimal bound on high stresses occurring between stiff fibers with arbitrary shaped cross sections, J. Math. Anal. Appl. 350 (2009), 306-312.




\bibitem{Yun-arXiv} K. Yun, An optimal estimate for electric fields on the shortest line segment between two spherical insulators in three dimensions, J. Differ. Equations 261 (2016), 148--188.





\end{thebibliography}
\end{document}